\documentclass[preprint,12pt]{elsarticle}

\usepackage{mathtools}
\usepackage{float}
\usepackage{tikz}
\usetikzlibrary{calc, math, patterns}
\usepackage{xcolor}
\usepackage{lineno}
\usepackage{amssymb}
\usepackage{amsmath}
\usepackage[hidelinks]{hyperref}
\usepackage[nameinlink]{cleveref}
\usepackage{calc}
\usepackage{mathrsfs}
\usepackage{booktabs}
\usepackage{bm}
\usepackage{subcaption}
\usepackage{siunitx}
\usepackage{pgfplots}
\usepackage{pgfplotstable}
\usepgfplotslibrary{groupplots,external}

\pgfplotsset{
    compat=1.18,
    filter discard warning=false,
    unbounded coords=discard}

\pgfplotsset{
    every axis title/.append style={yshift=-6pt},
    every x label/.append style={yshift=6pt},
    every y label/.append style={yshift=-6pt},
}

\definecolor{tabblue}{RGB}{31,119,180}
\definecolor{taborange}{RGB}{255,127,14}
\definecolor{tabgreen}{RGB}{44,160,44}
\definecolor{tabred}{RGB}{214,39,40}
\definecolor{tabpurple}{RGB}{148,103,189}
\definecolor{tabbrown}{RGB}{140,86,75}
\definecolor{tabpink}{RGB}{227,119,194}
\definecolor{tabgray}{RGB}{127,127,127}
\definecolor{tabolive}{RGB}{188,189,34}
\definecolor{tabcyan}{RGB}{23,190,207}

\pgfplotsset{
  tabblueplot/.style={
    color=tabblue,
    solid,
    mark=o,
    mark options={solid, draw=tabblue, fill=none}
  },
  taborangeplot/.style={
    color=taborange,
    solid,
    mark=+,
    mark options={solid, draw=taborange, fill=none}
  },
  tabgreenplot/.style={
    color=tabgreen,
    solid,
    mark=triangle,
    mark options={solid, draw=tabgreen, fill=none, rotate=180}
  },
  tabredplot/.style={
    color=tabred,
    solid,
    mark=diamond,
    mark options={solid, draw=tabred, fill=none}
  },
  tabpurpleplot/.style={
    color=tabpurple,
    solid,
    mark=triangle,
    mark options={solid, draw=tabpurple, fill=none, rotate=270}
  },
  tabbrownplot/.style={
    color=tabbrown,
    solid,
    mark=square,
    mark options={solid, draw=tabbrown, fill=none}
  },
  tabpinkplot/.style={
    color=tabpink,
    solid,
    mark=x,
    mark options={solid, draw=tabpink, fill=none}
  },
  tabgrayplot/.style={
    color=tabgray,
    solid,
    mark=pentagon,
    mark options={solid, draw=tabgray, fill=none}
  },
  taboliveplot/.style={
    color=tabolive,
    solid,
    mark=star,
    mark options={solid, draw=tabolive, fill=none}
  },
  tabcyanplot/.style={
    color=tabcyan,
    solid,
    mark=asterisk,
    mark options={solid, draw=tabcyan, fill=none}
  }
}

\pgfplotscreateplotcyclelist{mytab10}{
  tabblueplot\\
  taborangeplot\\
  tabgreenplot\\
  tabredplot\\
  tabpurpleplot\\
  tabbrownplot\\
  tabpinkplot\\
  tabgrayplot\\
  taboliveplot\\
  tabcyanplot\\
}
\usepackage{amsthm}
\newtheorem{remark}{Remark}

\crefname{section}{Section}{Sections}
\crefname{subsection}{Section}{Sections}
\crefname{subsubsection}{Section}{Sections}
\crefrangelabelformat{section}{#3#1#4--#5#2#6}
\crefrangelabelformat{subsection}{#3#1#4--#5#2#6}
\crefrangelabelformat{subsubsection}{#3#1#4--#5#2#6}

\newcommand{\grad}{\vec{\nabla}}
\newcommand{\tensor}[1]{\boldsymbol{#1}}

\newcommand{\matr}[1]{\bm{#1}}

\newcommand{\into}{\int_\Omega}

\newcommand{\matdisp}{\vec{u}}

\newcommand{\fludisp}{\vec{U}}
\newcommand{\seepage}{\vec{w}}

\newcommand{\volfrac}{\phi}
\newcommand{\fludensity}{\rho_f}
\newcommand{\soldensity}{\rho_s}
\newcommand{\biotcoeff}{\alpha}
\newcommand{\density}{\rho}
\newcommand{\aph}{(\alpha - \phi)}

\newcommand{\matstress}{\tensor{\sigma}''}

\newcommand{\splinespace}{\mathbb{S}}

\newcommand{\trialspaceu}{\mathscr{U}^h}
\newcommand{\trialspacep}{\mathscr{P}^h}
\newcommand{\trialspaceU}{\mathscr{V}^h}

\makeatletter
\nopreprintlinetrue
\makeatother

\begin{document}

\begin{frontmatter}

\title{Isogeometric Analysis for Explicit Wave Propagation in Poroelastic Media
}

\author[TUe]{Maarten M. Hodzelmans\corref{cor}}
\ead{m.m.hodzelmans@tue.nl}
\cortext[cor]{Corresponding author}
\author[TUe]{René R. Hiemstra}
\author[TUe]{Joris J.C. Remmers}
\author[TUe]{Clemens V. Verhoosel} 

\affiliation[TUe]{
    organization={Eindhoven University of Technology, Department of Mechanical Engineering},
    addressline={Traverse},
    postcode={PO Box 513, 5600 MB Eindhoven},
    country={The Netherlands}            
}

\begin{abstract}
For higher-order discretizations of explicit dynamics problems, Isogeometric Analysis (IGA) has several favorable properties as compared to classical Finite Element Analysis (FEA). While FEA produces spurious modes at orders beyond linear, this is not the case for IGA. Consequently, fewer degrees of freedom are required for comparable accuracy, larger timesteps can be taken, and the method is more robust for nonlinear problems. If outlier modes are removed, the timestep even becomes virtually independent of the order. In this paper, we investigate how these advantages apply to the poroelastic continuum model. \\
We consider both a primal formulation, wherein our variables are the displacement of the matrix material, the fluid displacement, and the pressure, as well as a reduced form wherein the pressure is eliminated. For our discretizations, we employ divergence-conforming spline spaces. Conforming spline spaces for the fluid displacement ensure inf-sup stability for the primal form, as well as a correct null space in the reduced form. Furthermore, we prove and demonstrate that the two formulations coincide when both displacements are discretized with conforming spline spaces. \\
Through spectral analysis, we find that the aforementioned benefits of IGA do carry over directly to the context of poroelasticity. In 1D, we split the discrete spectrum into fast and slow waves. When normalized against an analytical solution, each of these sub-spectra closely resembles results known in elasticity. Consequently, when poroelasticity is discretized with outlier-free IGA, the timestep is essentially independent of the order. We show this timestep scaling in 2D as well. \\

\end{abstract}

\begin{highlights}
\item The spectral advantages of IGA, as known in elasticity, carry over to poroelasticity
\item Outlier removal makes the critical timestep independent of the spline order
\item Divergence-conforming spline spaces result in a correct null space
\item With divergence-conforming spline spaces, the primal and reduced formulation coincide
\item One and two-dimensional poroelastodynamic benchmarks are studied
\end{highlights}


\end{frontmatter}

\section{Introduction}\label{chap:1}

Accurately modeling the propagation of seismic waves is of great importance in geophysics and subsurface engineering, for example in the context of hazard assessment \cite{edwardsSimulationsDevelopmentGround2019}. These models could help steer ongoing efforts to reinforce housing in regions that have become seismically active due to subsurface engineering activities, such as the Groningen region in the Netherlands \cite{muntendam-bosOverviewInducedSeismicity2022, vlekRiseReductionInduced2019}.
In addition, a model based on first principles is applicable to regions which are currently seismically inactive. This could help minimize the risks involved in future subsurface exploitation, such as geothermal energy production or CO$_2$ storage. In the context of risk assessment the shallow subsurface is of particular interest. When waves transition into the soft, uncompacted soils characteristic of the shallow subsurface, they increase in amplitude, an effect known as site amplification \cite{safakLocalSiteEffects2001}. Modeling this phenomenon is challenging, since the shallow subsurface exhibits strong material nonlinearity and heterogeneity.\\
Wave propagation in the shallow subsurface can be simulated through poroelastodynamics. This continuum description was proposed by Biot \cite{biotGeneralTheoryThreeDimensional1941, biotTheoryElasticityConsolidation1955}, building on Terzaghi's concept of an effective stress \cite{terzaghiBerechnugDurchlassigkeitTones1923, terzaghiErdbaumechanikAufBodenphysikalischer1925}. Biot's theory can capture the behavior of a porous matrix material saturated with a single-phase fluid, under the assumption that the microscopic constitution does not affect macroscopic behavior \cite{coussyPoromechanics2004}. It describes the spatial variation and temporal evolution of the matrix- and fluid displacements, $\vec{u}$ and $\vec{U}$ respectively, as well as the pore pressure $p$. To solve this model numerically, various formulations have been developed in the literature \cite{zienkiewiczDynamicBehaviourSaturated1984}. Typically, fluid inertia is neglected, naturally resulting in the displacement-pressure formulation ($u$-$p$). In this work, we instead study formulations that do include fluid inertia, since we intend to investigate its relevance for our application of near-surface seismic wave propagation. Therefore, we study both the full system, referred to as the $u$-$p$-$U$ formulation, as well as the reduced $u$-$U$ formulation. The former is applicable to both compressible and incompressible media, while the latter assumes compressibility \cite{zienkiewiczDynamicBehaviourSaturated1984}.\\
Isogeometric Analysis (IGA) is the concept of using the same geometric descriptions used in CAD for the purpose of analysis, effectively resulting in approximation of the physical fields using splines \cite{hughesIsogeometricAnalysisCAD2005, cottrellIsogeometricAnalysis2009}. These spline functions allow for higher-order continuity as compared to the Lagrange Finite Element Method (FEM), which has been shown to have various benefits. Notably, they are well-behaved at higher orders, as originally reported by Cottrel et al. \cite{cottrellIsogeometricAnalysisStructural2006, cottrellStudiesRefinementContinuity2007}. For elasticity problems, they investigated the discrete eigenfrequency spectrum. At polynomial orders beyond linear, Lagrange FEM results in both acoustic and optical eigenfrequencies, where only the former approximate physical modes. Meanwhile, spline bases at optimal regularity effectively eliminate such modes, instead computing only acoustic eigenfrequencies. Consequently, the number of Degrees of Freedom (DoFs) is reduced without reducing accuracy, which can be favorable for computational efficiency. Computing fewer eigenfrequencies also enables larger timesteps in an explicit solver. In addition to efficiency, the elimination of optical modes can improve robustness. As is remarked in Ref. \cite{evansExplicitHigherorderAccurate2018}, higher-order elements are impopular in explicit dynamics codes. Particularly in nonlinear simulations, where modes activate one another, accuracy across the entire spectrum is important. This is relevant when employing these methods to model softer soils.\\
Isogeometric analysis has been employed to solve static poroelasticity problems, using the mixed displacement pressure ($u$-$p$) formulation. It is well-established that this formulation can result in pressure oscillations if the timestep is too small \cite{vermeerAccuracyConditionConsolidation1981}. If IGA is used to enforce higher-order continuity in the pressure field, these oscillations are significantly reduced \cite{irzalIsogeometricFiniteElement2013}. The associated continuous pressure gradient enforces local conservation of mass as well.\\
Additional benefits are obtained by using a mixed formulation, wherein the pressure is discretized with lower-order functions as compared to those with which the displacement is discretized. This further reduces pressure oscillations resulting from too small timesteps, as well as pressure oscillations on account of a discontinuous permeability \cite{bekeleMixedMethodIsogeometric2022}. Similar conclusions with regards to pressure oscillations were drawn using isogeometric collocation with a mixed formulation \cite{morgantiMixedIsogeometricCollocation2018}.\\
Poroelasticity discretized through IGA has been studied in various other contexts. These include thermally coupled models \cite{bekeleIsogeometricAnalysisTHM2017}, the Darcy-Brinkman model \cite{vuongTwoFiniteElement2016}, and fracture mechanics \cite{irzalIsogeometricAnalysisBezier2014, fathiExtendedIsogeometricAnalysis2022, hagemanFlowNonNewtonianFluids2019, hagemanSubgridModelsMultiphase2020}. A comprehensive overview of advanced computational aspects in poroelasticity, including IGA, can be found in Ref. \cite{deborstComputationalMethodsFracture2017}.\\
The study of stable mixed formulations is important for the $u$-$p$-$U$ formulation as well. In the context of IGA, complexes of spline spaces have been constructed which conserve the structure of differential operators \cite{buffaIsogeometricDiscreteDifferential2011}. They thereby automatically satisfy the inf-sup condition, resulting in numerical stability \cite{evansIsogeometricDivergenceconformingBsplines2013}.\\
The application of IGA for dynamical problems can improve numerical efficiency when compared to FEM, and IGA has been employed to solve poroelasticity problems. However, to our best knowledge, the application of IGA to poroelastic wave propagation problems is unexplored. In this paper, we reduce this gap.\\
To this end, we perform a detailed comparison of the discrete spectra resulting from various Lagrangian FEM bases and B-spline bases. These spectra provide insight into the robustness of a method. In 1D, we segment eigenmodes into the fast and slow Biot waves, and normalize each cluster with an analytical solution. In 2D, we normalize against a refined reference. We also investigate how the choice of approximation space affects the critical timestep in an explicit time integration method. Finally, we study how and when the spectra resulting from the various poroelastodynamics formulations differ, and discuss structure-preserving discretizations. \\
Our paper is outlined as follows. In \autoref{chap:contiporo} we introduce the strong- and weak forms of the poroelasticity problem, which are discretized using IGA as outlined in \autoref{chap:discreteporo}. In \autoref{chap:1Dcase}, we discretize a 1D problem and analyze the resulting spectra, as well as the consequences for the critical timestep. This analysis is extended to 2D in \autoref{chap:2Dcase}. Our conclusions are summarized in \autoref{chap:conc}.

\section{Continuum poroelastodynamics}\label{chap:contiporo}
 
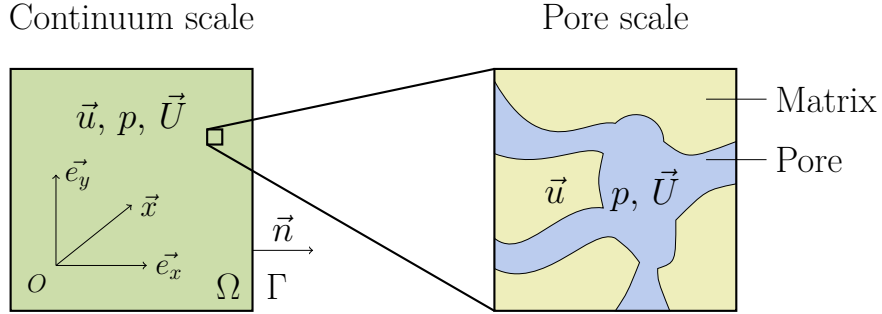
\begin{figure}[h!]
    \centering
    \makebox[0.8\linewidth][c]{%
        \scalebox{0.8}{%
            \tikzsetnextfilename{domainwithvars}
\begin{tikzpicture}

\useasboundingbox (-9, -3) rectangle (9, 3.5);

\def\ytopL{ 2}
\def\ybotL{-2}
\def\xleftL{-6}
\def\xrightL{-2}

\def\ytopR{ 2}
\def\ybotR{-2}
\def\xleftR{2}
\def\xrightR{6}

\def\offsetomega{0.1}

\definecolor{domaincol} {RGB}{204, 221, 170}
\definecolor{matrixcol} {RGB}{238, 238, 187}
\definecolor{channelcol}{RGB}{187, 204, 238}

\coordinate (LtopL) at (\xleftL,\ytopL);
\coordinate (LbotR) at (\xrightL,\ybotL);

\coordinate (RtopL) at (\xleftR,\ytopR);
\coordinate (RbotR) at (\xrightR,\ybotR);

\fill[domaincol] (LtopL) rectangle (LbotR);
\fill[matrixcol] (RtopL) rectangle (RbotR);

\draw[line width=1pt] (-2.75,  0.75) rectangle (-2.5, 1.0);
\draw[line width=1pt] (-2.75,  1.0) -- (RtopL);
\draw[line width=1pt] (-2.75,  0.75) -- (\xleftR, \ybotR);

\coordinate (AA) at (4.7, -1.3);
\coordinate (AB) at (4.7, -1.5);
\coordinate (AC) at (4.9, -2);
\coordinate (AD) at (4, -2);
\coordinate (AE) at (4.2, -1.7);
\coordinate (AF) at (4.4, -1.2);

\coordinate (BA) at (4.0, -0.8);
\coordinate (BB) at (3.6, -0.6);
\coordinate (BC) at (2.8, -1.7);
\coordinate (BD) at (2, -1.3);
\coordinate (BE) at (2, -0.8);
\coordinate (BF) at (2.5, -1.2);
\coordinate (BG) at (3.1, -0.2);
\coordinate (BH) at (3.8, -0.3);

\coordinate (CA) at (3.8, 0.6);
\coordinate (CB) at (3.5, 0.8);
\coordinate (CC) at (2.8, 0.3);
\coordinate (CD) at (2.0, 0.6);
\coordinate (CE) at (2, 1.8);
\coordinate (CF) at (2.8, 0.4);
\coordinate (CG) at (3.5, 1.2);
\coordinate (CH) at (4, 1.1);

\coordinate (DA) at (4.8, 0.8);
\coordinate (DB) at (5.2, 0.5);
\coordinate (DC) at (6, 0.8);
\coordinate (DD) at (6, 0.1);
\coordinate (DE) at (5.5, 0);
\coordinate (DF) at (5, -0.5);

\coordinate (EA) at (4.25, -1.0);
\coordinate (EB) at (3.7, 0.15);
\coordinate (EC) at (4.3, 1.4);
\coordinate (ED) at (4.8, 1.2);
\coordinate (EE) at (5.0, -1.2);

\path[fill=channelcol, draw=black]
  (AA) .. controls (AB) .. (AC) -- (AD)
       .. controls (AE) .. (AF)
       .. controls (EA) .. (BA)
       .. controls (BB) and (BC) .. (BD) -- (BE)
       .. controls (BF) and (BG) .. (BH)
       .. controls (EB) .. (CA)
       .. controls (CB) and (CC) .. (CD) -- (CE)
       .. controls (CF) and (CG) .. (CH)
       .. controls (EC) and (ED) .. (DA)
       .. controls (DB) .. (DC) -- (DD)
       .. controls (DE) .. (DF)
       .. controls (EE) .. (AA)
  -- cycle;

\draw[line width=1pt] (LtopL) rectangle (LbotR);
\draw[line width=1pt] (RtopL) rectangle (RbotR);

\coordinate (O) at ($(\xleftL,\ybotL) + (0.75,0.75)$);
\draw[->] (O) -- ($ (O) + (1.5, 0  ) $) node[right] {\large$\vec{e_x}$};
\draw[->] (O) -- ($ (O) + (0  , 1.5) $) node[right] {\large$\vec{e_y}$};
\draw[->] (O) -- ($ (O) + (1.25, 1) $) node[right] {\large$\vec{x  }$};
\node[below left] at (O) {$O$};

\node[above] at (\xleftL + 2, \ytopL+0.5) {\Large Continuum scale};
\node[above] at (\xleftR + 2, \ytopL+0.5) {\Large Pore scale};

\node at (\xrightL - 0.4, \ybotL + 0.4) {\Large $\Omega$};
\node at (\xrightL + 0.4, \ybotL + 0.4) {\Large $\Gamma$};
\draw[->] (\xrightL, -1) -- ++(1.0, 0) node[midway, above] {\Large $\vec{n}$};

\node at (\xleftL + 2.0, \ytopL - 0.75) {\Large  $\vec{u}, \, p,  \,  \vec{U} $};
\node at (3, 0) {\Large $\vec{u}$};
\node at (4.5, 0) {\Large $p, \, \vec{U}$};

\draw[] (\xrightR - 0.5, \ytopR - 0.5) -- ++(1.0, 0.0) node[right] {\Large Matrix};
\draw[] (\xrightR - 0.5, \ytopR - 1.5) -- ++(1.0, 0.0) node[right] {\Large Pore};

\end{tikzpicture}%
        }
    }
    \caption{Schematic of poromechanics. Pores are assumed to be tiny relative to the continuum scale, leading to three effective primary variables defined over the continuum domain $\Omega$.}
    \label{fig:domain}
\end{figure}

Consider a fluid-saturated porous medium within a domain $\Omega$, consisting on the pore-scale of a matrix and a fluid-filled pore space, modeled in a homogenized fashion, as illustrated in \autoref{fig:domain}. The domain has boundary $\Gamma$, with the outward-pointing normal vector $\vec{n}$. Over $\Omega$, we define the material point in the current configuration $\vec{x}$. The poroelastodynamic continuum model consists of three coupled equations for three independent variables: the displacement of the matrix material $\vec{u}(\vec{x})$, the pressure within the pores $p(\vec{x})$, and a pore-relative fluid displacement $\vec{w}(\vec{x})$. The latter is defined such that its time-derivative is the Darcy flux \cite{darcyFontainesPubliquesVille1856}. In addition, we define the absolute fluid displacement $\fludisp(\vec{x}) = \matdisp + \vec{w}/\volfrac$, with $\phi$ the fluid volume fraction. \\
We assume $\phi$, and all other material parameters, to be constant over $\Omega$, as heterogeneity is non-essential for our current study. Additionally, we assume that all variables are smooth over $\Omega$, and that displacements and displacement gradients are small.
We herein study the $u$-$p$-$U$ and $u$-$U$ formulations, which both allow for the inclusion of fluid inertia, in contrast to the commonly considered $u$-$p$ formulation \cite{zienkiewiczDynamicBehaviourSaturated1984}. In both cases, we consider a Lagrangian formulation. Note that either $\vec{U}$ or $\vec{w}$ can be chosen as a primary variable. While this choice is arbitrary in terms of physics, we prefer the former since it results in a sparser mass matrix, which is computationally beneficial. In the following, the two formulations are introduced, drawing inspiration from the derivation of Ref.~\cite{zienkiewiczDynamicBehaviourSaturated1984}. 

\subsection{Strong formulations}
The governing equations for poroelasticity, formulated in a Lagrangian manner, are the momentum balance in the homogenized material, the momentum balance in the fluid phase, and a resultant of mass conservation in the fluid phase, \emph{i.e.},
\begin{subequations}\label{eq:strongform_biot}
\begin{align}
    \grad \cdot \tensor{\sigma} + \density \vec{b} - \density \ddot{\matdisp} -\fludensity\left(\ddot{\seepage} + \dot{\seepage} \cdot \grad \dot{\seepage} \right) = 0,
    \\
    -\grad p + \fludensity \vec{b} = \fludensity \left(\ddot{\matdisp} + \frac{1}{\volfrac} \left(\ddot{\seepage} + \dot{\seepage} \cdot \grad \dot{\seepage} \right) \right) + \frac{\vec{R}}{\volfrac},
    \\
    \grad \cdot \dot{\seepage} + \biotcoeff \grad \cdot \dot{\matdisp} + \frac{1}{Q} \dot{p} = 0.
\end{align}
\end{subequations}
Here, $\tensor{\sigma}$ denotes the Cauchy stress in the homogenized material, $\vec{b}$ the body forces, and $\vec{R}$ the drag force due to seepage. The bulk density $\rho$ is a weighted average of the solid- and fluid density according to $\rho = \phi \fludensity + (1-\phi)\soldensity$. The parameter $\biotcoeff$ denotes the Biot-Willis coefficient \cite{biotElasticCoefficientsTheory1957}, and $Q$ the storage modulus.\\
To solve the system of equations \eqref{eq:strongform_biot}, two constitutive laws must be supplemented. The stress in the homogenized material is modeled as the Biot effective stress, meaning that it is decomposed as     
\begin{equation}\label{eq:biotstress}
    \tensor{\sigma}(\matdisp, p) = \matstress(\matdisp) - \biotcoeff p \boldsymbol{I},
\end{equation}
where $\matstress$ is the stress in the matrix material and $\tensor{I}$ the identity tensor. The first constitutive law is for the stress in the matrix, which we model through linear elasticity, \emph{i.e.},
\begin{equation}\label{eq:linela}
    \matstress = 2\mu\, \grad^s\matdisp + \lambda (\grad \cdot \matdisp) \, \boldsymbol{I},
\end{equation}
with $\lambda$ and $\mu$ the Lamé parameters and $\grad^s$ the symmetric gradient operator. The second constitutive law is Darcy's law \cite{darcyFontainesPubliquesVille1856}, with which resistance to seepage is modeled. In an isotropic setting, this is given as 
\begin{equation}\label{eq:darcy}
    \vec{R} = \frac{\volfrac}{k}\dot{\vec{w}},
\end{equation}
with $k$ the Darcy conductivity, in the terminology proposed in Ref.~\cite{kumpelPoroelasticityParametersReviewed1991}. This conductivity is defined as the ratio of the permeability and the dynamic viscosity. Note, however, that these never occur independently within poroelastodynamic theory.\\    
To express the strong form for the $u$-$p$-$U$ formulation, we neglect the convective terms $\dot{\seepage} \cdot \grad \dot{\seepage}$, since, as Zienkiewicz and Shiomi \cite{zienkiewiczDynamicBehaviourSaturated1984} note, in the setting of earthquake engineering this term is typically smaller than the uncertainty with which the parameters are defined. Furthermore, we neglect body forces. Upon substitution of Darcy's law \eqref{eq:darcy}, the definition of the Biot effective stress \eqref{eq:biotstress}, and the definition of $\fludisp$ into \eqref{eq:strongform_biot}, after rearrangement, the $u$-$p$-$U$ formulation is obtained as
\begin{subequations}\label{eq:threefield_strong}
\begin{align}\label{eq:threefield_strong_1}
    \grad \cdot \matstress(\matdisp) - (\biotcoeff - \volfrac) \grad p - (1 - \volfrac) \soldensity \ddot{\matdisp} + \frac{\volfrac^2}{k}\left(\dot{\fludisp} - \dot{\matdisp}\right) = 0,
    \\\label{eq:threefield_strong_2}
    \volfrac \grad p + \volfrac \fludensity \ddot{\fludisp} + \frac{\volfrac^2}{k}\left(\dot{\fludisp} - \dot{\matdisp}\right) = 0,
    \\\label{eq:threefield_strong_3}
    \volfrac \grad \cdot \fludisp + (\biotcoeff - \volfrac) \grad \cdot \matdisp + \frac{1}{Q} p = p_0(\vec{x}),
\end{align}
\end{subequations}
where the initial pressure $p_0(\vec{x})$ appears as an integration constant. For the sake of brevity, we assume this vanishes from here on. Note that the pressure contributes to the momentum balance with a factor $\alpha$, through the Biot effective stress \eqref{eq:biotstress}, and additionally with a factor $\volfrac$, as a consequence of volume-averaging the stresses in the constituents. These are combined in the term with $(\biotcoeff - \volfrac)$ in \eqref{eq:threefield_strong_1}.\\
The $u$-$U$ formulation follows directly from here. If the homogenized material is compressible (which implies that $Q$ is finite), we can rearrange the mass balance \eqref{eq:threefield_strong_3} explicitly for the pressure, as 
\begin{equation}\label{eq:p_elim_analytical}
    p = -Q \left(\volfrac \grad \cdot \fludisp + (\biotcoeff - \volfrac) \grad \cdot \matdisp \right).
\end{equation}
The pressure is then eliminated through substitution into \eqref{eq:threefield_strong_1} and \eqref{eq:threefield_strong_2}, resulting in the $u$-$U$ formulation:
\begin{subequations}\label{eq:twofield_strong}
\begin{align}
\begin{split}
        \grad \cdot \matstress(\matdisp) - (1 - \volfrac) \soldensity \ddot{\matdisp} \, \, + \hphantom{ 0}&
        \\
        (\biotcoeff - \volfrac) Q\grad \left(\volfrac \grad \cdot \fludisp + (\biotcoeff - \volfrac) \grad \cdot \matdisp \right) + \frac{\volfrac^2}{k}\left(\dot{\fludisp} - \dot{\matdisp}\right) = 0, &
    \end{split}
    \\
    \begin{split}
        \volfrac \fludensity \ddot{\fludisp} -         
        \volfrac Q \grad  \left(\volfrac \grad \cdot \fludisp + (\biotcoeff - \volfrac) \grad \cdot \matdisp \right) + \frac{\volfrac^2}{k}\left(\dot{\fludisp} - \dot{\matdisp}\right) = 0. &
    \end{split}
\end{align}
\end{subequations}
Both systems \eqref{eq:threefield_strong} and \eqref{eq:twofield_strong} are complemented by initial conditions, $\matdisp(\vec{x}, 0) = f(\vec{x}), \, \fludisp(\vec{x}, 0) = g(\vec{x})$, and by two pairs of boundary conditions. Recall that an initial pressure $p_0$ may be defined as well, which we neglect here. At every point of the boundary, either the displacement or the traction $\vec{t} = \sigma \cdot \vec{n}$ must be specified, as well as either the fluid displacement or the pore pressure, \emph{i.e.}
\begin{subequations}\label{eq:BCs}
\begin{align}
    \matdisp = \matdisp^D \quad \mathrm{on} \quad \Gamma^D_s, \qquad 
    \vec{t} = \vec{t}^{\,\,N} \quad \mathrm{on} \quad \Gamma^N_s,
    \\
    \fludisp = \fludisp^D \quad \mathrm{on} \quad \Gamma^D_f, \qquad 
    p = p^N \quad \mathrm{on} \quad \Gamma^N_f.
\end{align}
\end{subequations}
Herein, the boundary is split into the non-overlapping Dirichlet and Neumann boundary for both pairs. That is, 
$\Gamma = \Gamma^D_s \, \cup \Gamma_s^N = \Gamma^D_f \cup \Gamma^N_f$ 
and 
$\Gamma^D_s \, \cap \Gamma_s^N = \Gamma^D_f \cap \Gamma^N_f = \emptyset$. In multi-D, the boundary may be split into Dirichlet and Neumann conditions differently for each vectorial component. For the sake of notational brevity, this is not done here.
We explicitly note that the traction is defined through the Biot effective stress \eqref{eq:biotstress}, rather than exclusively the stress in the matrix.

\subsection{Weak formulations}

Weak forms of both the $u$-$p$-$U$-formulation \eqref{eq:threefield_strong} and the $u$-$U$-formulation\eqref{eq:twofield_strong} are obtained through the method of weighted residuals \cite{zienkiewiczChapter3Weak2013}. For that purpose, we define three test functions, $\vec{v}(\vec{x}),\, \vec{\eta}(\vec{x})$ and $\psi(\vec{x})$. The first two test functions correspond to the solid- and fluid momentum equations respectively, and the last to the pressure equation. The balance equations are multiplied by their corresponding test function, and the resulting statement is required to hold for any suitable test function. Subsequently integrating over the domain, integrating by parts, and applying Gauss' theorem yields the weak formulations. In terms of linear- and bilinear operators, the weak $u$-$p$-$U$ formulation is written as
\begin{subequations}\label{eq:threefield_weak}
\begin{align}
    \mathrm{Find }\,(\vec{u},\vec{p},\vec{U}) \in
    (\mathscr{U}\times\mathscr{P}\times\mathscr{V})
    \,\,\mathrm{such\ that:}
    \qquad\qquad\qquad
    \notag\\
\begin{split}
    \mathcal{K}(\vec{v}, \matdisp)
    -\mathcal{G}(\vec{v},p)
    +\mathcal{M}_s(\vec{v},\ddot{\matdisp})
    -\mathcal{C}(\vec{v},\dot{\fludisp}-\dot{\matdisp})
    =
    \\
    \mathcal{F}_t(\vec{v})+\mathcal{F}_p(\vec{v})
    \quad&\forall\,\vec{v}\in\mathscr{W},
\end{split}
    \label{eq:threefield_weak_1}
    \\
    -\mathcal{H}(\vec{\eta},p)
    +\mathcal{M}_f(\vec{\eta},\ddot{\fludisp})
    +\mathcal{C}(\vec{\eta},\dot{\fludisp}-\dot{\matdisp})
    =
    -\mathcal{F}_p(\vec{\eta})
    \quad&\forall\,\vec{\eta}\in\mathscr{E},
    \label{eq:threefield_weak_2}
    \\
    \mathcal{H}(\fludisp,\psi)
    +\mathcal{G}(\matdisp,\psi)
    +\mathcal{P}(\psi,p)
    =0
    \quad&\forall\,\psi\in\mathscr{S}.
    \label{eq:threefield_weak_3}
\end{align}
\end{subequations}
Herein, $\mathscr{U}, \mathscr{P}$ and $\mathscr{V}$, are the trial spaces and $\mathscr{W}, \mathscr{E}$, and $\mathscr{S}$ the test spaces. The individual operators in this weak form are defined as
\begin{subequations}\label{eq:operators1}
\begin{align}
    \mathcal{K}(\vec{v}, \matdisp) 
    &:= \lambda
    \into 
    \grad \cdot \vec{v} \, \grad \cdot \matdisp 
    \, dV 
    + 
    2\mu 
    \into 
    \grad \vec{w} : (\grad^s\matdisp)
    \, dV,
    \\
    \mathcal{M}_s(\vec{v}, \ddot{\matdisp}) 
    &:= (1 - \phi) \rho_s 
    \into
    \vec{v} \cdot \ddot{\matdisp}
    \, dV,
    \\
    \mathcal{M}_f(\vec{\eta}, \ddot{\fludisp}) 
    &:= \volfrac \rho_f
    \into
    \vec{\eta} \cdot \ddot{\fludisp}
    \, dV,
    \\
    \mathcal{C}(\vec{v}, \dot{\matdisp})
    &:= \frac{\volfrac^2}{k}
    \into 
    \vec{v} \cdot \dot{\matdisp}
    \, dV,
    \label{eq:dampingoperator}
    \\
    \mathcal{G}(\vec{v}, p) 
    &:= (\biotcoeff - \volfrac) 
    \into
    \grad \cdot \vec{v} \, p
    \, dV,
    \\
    \mathcal{H}(\vec{\eta}, p)
    &:= \volfrac
    \into 
    \grad \cdot \vec{\eta} \, p
    \, dV,
    \\
    \mathcal{P}(\psi, p)
    &:= \frac{1}{Q}
    \into
    \psi \, p\
    \, dV,
    \\
    \mathcal{F}_t(\vec{v})
    &:= \int_{\Gamma^N}
    \vec{v} \cdot \vec{t}
    \, dS,
    \label{eq:tractBC}
    \\
    \mathcal{F}_p(\vec{v})
    &:= \volfrac \int_{\Gamma^N}
    \vec{v} \cdot p\vec{n} 
    \, dS .
    \label{eq:pBC}
\end{align}  
\end{subequations}
Suitable test- and trial functions are found in those Sobolev spaces $\mathbf{H}$ for which the operators \eqref{eq:operators1} remain finite. In addition, the Dirichlet boundary conditions \eqref{eq:BCs} are incorporated directly into the trial space, and by requiring that test functions vanish on the Dirichlet boundary. Hence, we define these spaces as
\begin{subequations}\label{eq:spaces_continuous}
\begin{align}
    \mathscr{U} = 
    \{\vec{u} \,|\, \vec{u}\in \mathbf{H}^1(\Omega), \,
    \vec{u} = \vec{u}^D \mathrm{on}\, \Gamma^D_s \},
    \\
    \mathscr{W} = 
    \{\vec{v} \,|\, \vec{v}\in \mathbf{H}^1(\Omega), \,
    \vec{v} = \vec{0}\, \mathrm{on}\, \Gamma^D_s \},
    \\
    \mathscr{V} = 
    \{\vec{U} \,|\, \vec{U}\in \mathbf{H}(\mathrm{div},\Omega), \,
    \vec{U} = \vec{U}^D\, \mathrm{on}\, \Gamma^D_f \},
    \\
    \mathscr{E} = 
    \{\vec{\eta} \,|\, \vec{\eta}\in \mathbf{H}(\mathrm{div},\Omega), \,
    \vec{\eta} = \vec{0} \, \mathrm{on}\, \Gamma^D_f \}.
    \\
    \mathscr{P} = 
    \{p \,|\, p\in \mathbf{L}^2(\Omega) \},    
    \\
    \mathscr{S} = 
    \{\psi \,|\, \psi\in \mathbf{L}^2(\Omega)\},
\end{align}
\end{subequations}
The weak form of the $u$-$U$ formulation is obtained through the same procedure, resulting in 
\begin{subequations}\label{eq:twofield_weak}
\begin{align}    
    \mathrm{Find } \, (\vec{u},\vec{U}) \in (\mathscr{U} \times \mathscr{W})\ 
    \mathrm{such}\ 
    \mathrm{that:}
    \qquad\qquad\qquad\qquad\qquad
    \notag    \\
\begin{split}
    \mathcal{K}(\vec{v}, \vec{u}) + 
    \mathcal{R}_1(\vec{v}, \matdisp) + 
    \mathcal{R}_2(\vec{v}, \fludisp) + 
    \mathcal{M}_s(\vec{v}, \ddot{\matdisp}) -
    \mathcal{C}(\vec{v}, \dot{\fludisp} - \dot{\matdisp}) 
    = 
    \\
    \mathcal{F}_t(\vec{v}) + 
    \mathcal{F}_p(\vec{v})
    \quad &\forall \, \vec{v} \in \mathscr{W},
\end{split}
\label{eq:twofield_weak_1}
\\
\mathcal{R}_3(\vec{\eta}, \fludisp) + 
\mathcal{R}_2(\vec{\eta}, \matdisp) + 
\mathcal{M}_f(\vec{\eta}, \ddot{\fludisp}) + 
\mathcal{C}(\vec{\eta}, \dot{\fludisp} - 
\dot{\matdisp}) 
=
- \mathcal{F}_p(\vec{\eta}) 
\quad &\forall \, \vec{\eta} \in \mathscr{E}.
\label{eq:twofield_weak_2}
\end{align}
\end{subequations}
Most operators were already defined for the $u$-$p$-$U$ system in \eqref{eq:operators1}. The test- and trial spaces are those defined in \eqref{eq:spaces_continuous}. The distinction is that rather than coupling to the pressure, we have additional coupling between the two displacements through the operators $\mathcal{R}$, defined as    
\begin{subequations}\label{eq:operators2}
\begin{align}
    \mathcal{R}_1(\vec{v}, \matdisp) 
    :=& (\biotcoeff - \volfrac)^2 Q
    \into
    \grad \cdot \vec{v} \,\, \grad \cdot \matdisp
    \, dV,
    \\
    \mathcal{R}_2(\vec{v}, \fludisp) 
    :=& (\biotcoeff - \volfrac) \volfrac \, Q
    \into
    \grad \cdot \vec{v} \,\, \grad \cdot \fludisp
    \, dV,
    \\
    \mathcal{R}_3(\vec{\eta}, \fludisp) 
    :=& 
    \volfrac^2 Q
    \into
    \grad \cdot \vec{\eta} \,\, \grad \cdot \fludisp
    \, dV.
\end{align}
\end{subequations}
We refer the reader to Ref. \cite{phillipsCouplingMixedContinuous2007} for a detailed exposition of the functional setting, including well-posedness proofs.

\section{Discrete poroelastodynamics}\label{chap:discreteporo}

We discretize the formulations of \autoref{chap:contiporo} in space through the Bubnov-Galerkin procedure \cite{zienkiewiczChapter3Weak2013}. We introduce three sets of basis functions $\matr{N}^u$, $\matr{N}^U$ and $\matr{N}^p$, where the former two are vector-valued. Their spans form our discrete trial spaces, which are a subset of their continuous counterparts,
\begin{align}
    \mathscr{U}^h = \mathrm{span}\{\matr{N}^u\}, \quad
    \mathscr{V}^h = \mathrm{span}\{\matr{N}^U\}, \quad
    \mathscr{P}^h = \mathrm{span}\{\matr{N}^p\}.
\end{align}
We then search for a solution within these trial spaces, yielding the approximations
\begin{subequations}\label{eq:galerkinapprox}
\begin{align}
    \matdisp \approx \matdisp^h = 
    \sum_{I=1}^n \vec{N}_I^{\,u} u_I 
    &= [\vec{N}^{\,u}_1,\, ...,\, \vec{N}^{\,u}_n] \matr{u}
    \\
    \fludisp \approx \fludisp^h =
    \sum_I \vec{N}^{\,U}_I U_I 
    &= [\vec{N}^{\,U}_1,\, ...,\, \vec{N}^{\,U}_n] \matr{U}
    \\
    p \approx p^h = \sum_I N^{\,p}_I p_I 
    &= [N^{\,p}_1,\, ...,\, N^{\,p}_n] \matr{p}
\end{align}
\end{subequations}
where the arrays $\matr{u}$, $\matr{U}$, and $\matr{p}$ contain DOFs for the solution. The test functions are similarly discretized. \\
Note that we do not explicitly write boundary conditions in the discrete setting for notational brevity. Nonetheless, we only search for solutions in the subspaces of $\mathscr{U}^h, \mathscr{P}^h$ and $\mathscr{V}^h$ which satisfy the boundary conditions, and only test against the subspaces which vanish on the Dirichlet boundaries.\\
Both formulations are discretized in space by substituting the basis functions into the operators \eqref{eq:operators1} and \eqref{eq:operators2}. For the $u$-$p$-$U$ formulation, this yields the coupled system
\begin{align}\label{eq:upUfullsys}
\begin{split}
    \begin{bmatrix}
        \matr{M}_s & 0 & 0\\
        0 & 0 & 0\\
        0 & 0 & \matr{M}_f
    \end{bmatrix}
    \begin{bmatrix}
        \ddot{\matr{u}} \\ \ddot{\matr{p}} \\ \ddot{\matr{U}}
    \end{bmatrix} 
    +    
    \begin{bmatrix}
        \matr{C}_1 & 0 & -\matr{C}_2 \\
        0 & 0 & 0\\
        -\matr{C}_2^T & 0 & \matr{C}_3
    \end{bmatrix}
    \begin{bmatrix}
        \dot{\matr{u}} \\ \dot{\matr{p}} \\ \dot{\matr{U}}
    \end{bmatrix}
    + \\
    \begin{bmatrix}
        \matr{K} & -\matr{G} & 0\\
        -\matr{G}^T & \matr{P} & -\matr{H}^T\\
        0 & -\matr{H} & 0
    \end{bmatrix}
    \begin{bmatrix}
        \matr{u} \\ \matr{p} \\ \matr{U}
    \end{bmatrix}
    =
    \begin{bmatrix}
        \matr{f}_t + \matr{f}_s
        \\
        0
        \\
        -\matr{f}_f
    \end{bmatrix}
    ,
\end{split}    
\end{align}    
using the discretized operators
\begin{subequations}\label{eq:discreteoperators}
\begin{align}
    [\matr{M}_s]_{IJ} &= \mathcal{M}_s(\vec{N}_I^{\,u}, \vec{N}_J^{\,u}), 
    &
    [\matr{K}  ]_{IJ} &= \mathcal{K}  (\vec{N}_I^{\,u}, \vec{N}_J^{\,u}), 
    \\
    [\matr{M}_f]_{IJ} &= \mathcal{M}_f(\vec{N}_I^{\,U}, \vec{N}_J^{\,U}),
    &
    [\matr{P}  ]_{IJ} &= \mathcal{K}  (N_I^p, N_J^p), 
    \\
    [\matr{G}  ]_{IJ} &= \mathcal{G  }(\vec{N}_I^{\,u}, N_J^p), 
    &
    [\matr{H}  ]_{IJ} &= \mathcal{H}  (\vec{N}_I^{\,u}, N_J^p), 
    \\
    [\matr{C}_1]_{IJ} &= \mathcal{C}  (\vec{N}_I^{\,u}, \vec{N}_J^{\,u}), 
    &
    [\matr{R}_1]_{IJ} &= \mathcal{R}_1(\vec{N}_I^{\,u}, \vec{N}_J^{\,u}), 
    \\
    [\matr{C}_2]_{IJ} &= \mathcal{C}  (\vec{N}_I^{\,U}, \vec{N}_J^{\,u}), 
    &
    [\matr{R}_2]_{IJ} &= \mathcal{R}_2(\vec{N}_I^{\,U}, \vec{N}_J^{\,u}), 
    \\
    [\matr{C}_3]_{IJ} &= \mathcal{C}  (\vec{N}_I^{\,U}, \vec{N}_J^{\,U}), 
    &
    [\matr{R}_3]_{IJ} &= \mathcal{R}_3(\vec{N}_I^{\,U}, \vec{N}_J^{\,U}), 
    \\
    [\matr{f}_s]_{I } &= \mathcal{F}_p(\vec{N}_I^{\,u}), 
    &
    [\matr{f}_t]_{I } &= \mathcal{F}_t(\vec{N}_I^{\,u}),
    \\
    [\matr{f}_f]_{I } &= \mathcal{F}_p(\vec{N}_I^{\,U}).
\end{align}
\end{subequations}
We note immediately that the discrete version of the mass balance \eqref{eq:threefield_weak_3}, which forms the middle row of \eqref{eq:upUfullsys}, may be written explicitly for the pressure,
\begin{equation}\label{eq:discrmassbalance}
    \matr{p} = -\matr{P}^{-1}\left(\matr{G}^T \matr{u} + \matr{H}^T \matr{U}\right).
\end{equation}
This is possible so long as $\matr{P}$ is invertible. That is the case if the storage modulus $Q$ is finite, which mirrors the condition we had for reducing to the $u$-$U$ formulation algebraically. With \eqref{eq:discrmassbalance}, the pressure is condensed out of \eqref{eq:upUfullsys}. This results in a form that closely resembles the discrete $u$-$U$ formulation. Anticipating this, after condensation we write
\begin{equation}\label{eq:secondorderform}
\underbrace{
    \begin{bmatrix}
        \matr{M}_s & \matr{0} \\
        \matr{0}   & \matr{M}_f
    \end{bmatrix}
}_{:= \overline{\matr{M}}}
    \begin{bmatrix}
        \ddot{\matr{u}} \\
        \ddot{\matr{U}}
    \end{bmatrix}
     + 
\underbrace{
    \begin{bmatrix}
        \matr{C}_1 & -\matr{C}_2 \\
        -\matr{C}_2^T   & \matr{C}_3
    \end{bmatrix}
}_{:= \overline{\matr{C}}}
    \begin{bmatrix}
        \dot{\matr{u}} \\
        \dot{\matr{U}}
    \end{bmatrix}
     +         
    \overline{\matr{K}}
    \begin{bmatrix}
        {\matr{u}} \\
        {\matr{U}}
    \end{bmatrix}
     = 
     \begin{bmatrix}
         \matr{f}_t + \matr{f}_s\\
         -\matr{f}_f
     \end{bmatrix},
\end{equation}
where the system stiffness matrix $\overline{\matr{K}}$ for the $u$-$p$-$U$ formulation is
\begin{equation}\label{eq:upUstiff}
    \overline{\matr{K}}_{u\text{-}p\text{-}U} = 
    \begin{bmatrix}
      \matr{K}  + \matr{G} \matr{P}^{-1} \matr{G}^T  
            &              \matr{G} \matr{P}^{-1} \matr{H}^T \\
      \qquad      \matr{H} \matr{P}^{-1} \matr{G}^T  
            &              \matr{H} \matr{P}^{-1} \matr{H}^T
    \end{bmatrix}.
\end{equation}
Note that the reduced mass matrix $\overline{\matr{M}}$ is positive definite, while the mass matrix before condensation in \autoref{eq:upUfullsys} is singular due to the zero-block on the diagonal. Thus, the generalized eigenvalue problem corresponding to \autoref{eq:upUfullsys} is singular, while that corresponding to the condensed form \autoref{eq:secondorderform} is non-singular. For this reason, in the remainder of this work, we only consider the $u$-$p$-$U$ formulation in this condensed form.\\
The discrete $u$-$U$ system consists of the same mass, damping and forcing terms \eqref{eq:secondorderform}. The only difference is in the system stiffness. For the $u$-$U$ formulation, this is
\begin{equation}\label{eq:uUstiff}
    \overline{\matr{K}}_{u\text{-}U} = 
    \begin{bmatrix}
        \matr{K} + \matr{R}_1 
        &          \matr{R}_2  \\
        \qquad     \matr{R}_2^T
        &          \matr{R}_3
    \end{bmatrix}.
\end{equation}
Note that the stiffness matrices too have a shared structure. This similarity stems from the fact that the exact same physics and assumptions underlie the systems. The distinction is in how we embed the mass balance into the two momentum balances. This is either done algebraically in the strong form \eqref{eq:p_elim_analytical}, or through linear algebra in the discrete setting \eqref{eq:discrmassbalance}. In the latter case the pressure is explicitly discretized, rather than inherited from the displacement discretizations. This ensures that a compatible discretization can be chosen. It does, however, introduce additional computational effort, since the matrix $\matr{P}$ must be inverted. Conversely, the discrete $u$-$U$ formulation does not require inverting this matrix, but does not allow for a choice in pressure discretization either.

\begin{remark}
    For the analysis presented in this work, we consider only the case where the bulk is compressible. If that is not the case, neither of the forms presented in the above are valid: both $\overline{\matr{K}}_{u\text{-}p\text{-}U}$ and $\overline{\matr{K}}_{u\text{-}U}$ will approach infinity. While the non-reduced $u$-$p$-$U$ formulation is applicable to incompressible media, this is outside the scope of this work.
\end{remark}   

\subsection{B-spline basis functions}
Starting in the univariate setting, we choose B-splines for our bases $\matr{N}$ (\autoref{eq:galerkinapprox}), which are piece-wise polynomial functions that allow for higher-order regularity \cite{cottrellNURBSPreAnalysisTool2009}. We denote a spline space as 
\begin{equation}
    \mathbb{S}_k^p = \mathrm{span}\left\{N_{i, p}(x) \right\}_{i=1}^{n},
\end{equation}
where $p$ denotes the polynomial order and $k$ the number of times the space is continuously differentiable. To construct such a space, a knot vector $\Xi$ is defined, containing non-decreasing real-valued knot values $\xi_i$. The basis is constructed in the parameter domain with the coordinate $\zeta$, and subsequently mapped onto $\Omega$. The $i$-th spline $N_i$ can be evaluated using the recursive Cox- de Boor formulas \cite{coxNumericalEvaluationBsplines1972, deboorCalculatingBsplines1972}. That is, piece-wise constant basis functions are defined as
\begin{equation}
    N_{i,0}(\zeta) =
    \begin{cases}
    1, & \text{if } \xi_i \le \zeta < \xi_{i+1}, \\
    0, & \text{otherwise.}
    \end{cases}
\end{equation}
Higher-order basis functions are defined recursively as
\begin{equation}
    N_{i, p}(\zeta) = \frac{\zeta-\xi_i}{\xi_{i+p} - \xi_i}N_{i, p-1}(\zeta) + \frac{\xi_{i+p+1}-\zeta}{\xi_{i+p+1} - \xi_{i+1}}N_{i+1, p-1}(\zeta).
\end{equation}
We consider open knot vectors, meaning that the first- and last values are repeated $p+1$ times, resulting in interpolating functions at the domain boundaries. Such a basis has continuity $k=p-1$ (optimal regularity) if none of the interior values in the knot vector are repeated. Conversely, repeating an interior value once introduces one DOF and reduces the continuity of the basis over that knot by one. Thus, by repeating all interior knots the continuity of the entire basis is reduced. At continuity zero, a basis is obtained that spans the same space as a Lagrange basis.\\
In the following, we present a 2D problem. To discretize it, a multivariate B-spline space in 2D is defined simply as the tensor product of two univariate spaces, \emph{i.e.},
\begin{equation}
    \mathbb{S}^{p, p}_{k, k} = \mathbb{S}_k^p \otimes \mathbb{S}_k^p.
\end{equation}
IGA as defined in Ref. \cite{hughesIsogeometricAnalysisCAD2005} uses Non-Uniform Regular B-Splines (NURBS) for analysis, since this superset of B-splines allows exact descriptions of conic sections. However, as this is non-essential to our study, we consider only B-splines.

\subsection{Optimal spline spaces}\label{sec:outliers}
When using optimal regularity B-splines for analysis, order-elevation improves the accuracy of the entire discrete eigenvalue-spectrum, with one exception. When the polynomial degree is two or more, outlier modes appear \cite{cottrellIsogeometricAnalysisStructural2006}. These modes are non-physical, and with a uniformly spaced knot vector, will have frequencies significantly larger than the remainder of the spectrum. Consequently, they will constrain the maximum stable timestep in an explicit solver significantly.\\
The existence of such outlier modes has been well-established in the literature \cite{cottrellIsogeometricAnalysisStructural2006}, and will be demonstrated in the remainder of this work as well. These modes originate from the boundaries of the domain. In a 1D elasticity problem, the characteristic solutions have vanishing even derivatives on a Dirichlet boundary, and vanishing odd derivatives on a Neumann boundary. Meanwhile, a higher-order spline basis does admit such derivatives \cite{hiemstraRemovalSpuriousOutlier2021}. As these are not natural solutions to the problem, their appearance in any solution is "penalized" by high eigenvalues. A reduced basis, wherein these derivatives are eliminated, results in an outlier-free discrete system \cite{hiemstraRemovalSpuriousOutlier2021}. Such bases, known as optimal spline spaces, were first studied by Takacs and Takacs \cite{takacsApproximationErrorEstimates2016}, and their optimal approximation properties were proven in, among others, Refs.~\cite{floaterOptimalSplineSpaces2017, floaterOptimalSplineSpaces2019}.\\
There are multiple approaches to numerically perform this outlier elimination, of which an overview is provided in Ref. \cite{voetMassLumpingOutlier2026}. The procedure we use here is that of Hiemstra et al. \cite{hiemstraRemovalSpuriousOutlier2021}. First, two submatrices for Dirichlet and Neumann boundaries are constructed,
\begin{subequations}
\begin{align}
    [\matr{A}_D]_{ij} = \Delta_n^{i} \left(\left.N_j\right|_{\Gamma_D} \right)
    \qquad & i \in \{ 1, 2, 3, ..., (p-1)/2\},
    \\
    [\matr{A}_N]_{ij} = \Delta_n^{(2i-1)/2} \left(\left.N_j\right|_{\Gamma_N} \right)
    \qquad & i \in \{1, 2, 3, ..., p/2\},
\end{align}
\end{subequations}
where $\Delta_n$ denotes the second-order derivative with respect to the normal of the boundary. Each row in these matrices, upon multiplication by a solution vector, evaluates a normal derivative at the boundary that we intend to constrain to zero. To that end, for both matrices a nullspace $\matr{C}$ is constructed. In the first appendix of Ref. \cite{hiemstraRemovalSpuriousOutlier2021}, an algorithm is provided that does so while retaining non-negativity, minimal compact support and partition of unity. The two nullspaces together form the constraint matrix 
\begin{equation}
\matr{C} = 
    \begin{bmatrix}
        \matr{C}_D & 0        & 0            \\
        0          & \matr{I} & 0            \\
        0          & 0        & \matr{C}_N
    \end{bmatrix}.
\end{equation}
This projects the basis into a lower-dimensional space wherein the constraints are satisfied. That is, the reduced basis is found as $\matr{N}^- = \matr{C} \matr{N}$. In practice, the constraint matrix is computed once, before runtime, and it is applied through left- and right multiplication with the linear system during the solution routine, making it non-intrusive.

\subsection{Conforming spline spaces}\label{subsec:conforming}
For the $u$-$p$-$U$ formulation we need to address compatibility, as not all combinations of $\matr{N}^u$, $\matr{N}^p$ and $\matr{N}^U$ result in a stable and converging method. While we consider a rigorous analysis of stability and optimality beyond the scope of the current work, we do briefly discuss what compatible bases might look like. Additionally, in \autoref{chap:1Dcase}, we demonstrate the consequences of choosing incompatible bases. \\
In the general case, a structure preserving discretization can be derived on the basis of finite element exterior calculus \cite{arnoldFiniteElementExterior2018, hiemstraHighOrderGeometric2014}. Specific to the poroelasticity problem, requirements on compatible discretizations are noted in Ref. \cite{phillipsCouplingMixedContinuous2007}. For our $u$-$p$-$U$ weak form \eqref{eq:threefield_weak}, the key requirement is compatibility between the discrete fluid displacement space and the pressure space. The divergence of the former should be contained in the latter. \\
A suitable class of IGA-elements are the divergence-conforming ones from the sequence proposed in Ref. \cite{buffaIsogeometricDiscreteDifferential2011} (see also their application in, among others, Refs.~\cite{evansIsogeometricDivergenceconformingBsplines2013, evans_divconforming_darcystokesbrinkman_2013}). Accordingly, in 2D, we choose our solution spaces $\mathscr{V}^h$ and $\mathscr{P}^h$ for $\fludisp$ and $p$ as
\begin{subequations}\label{eq:conforming2D}
\begin{align}
    \trialspaceU &= 
    \mathbb{S}^{k_1, k_2-1}_{\alpha_1, \alpha_2-1} 
    \times 
    \mathbb{S}^{k_1-1, k_2}_{\alpha_1-1, \alpha_2},
    \label{eq:conforming2D-vector}
    \\
    \trialspacep &=
    \mathbb{S}^{k_1-1, k_2-1}_{\alpha_1-1, \alpha_2-1},
    \label{eq:conforming2D-scalar}
\end{align}
\end{subequations}
where $\alpha_i \geq 0$ and $k_i > \alpha_i$. It is easily seen that the divergence of the former space equals the second space. In 1D, these spaces reduce to
\begin{subequations}\label{eq:conforming1D}
\begin{align}
    \trialspaceU =& 
    \mathbb{S}^{k}_\alpha,
    \\
    \trialspacep =&
    \mathbb{S}^{k-1}_{\alpha-1}.
\end{align}
\end{subequations}
Beyond compatibility, we prove in \ref{app:equivalence} that the $u$-$p$-$U$ formulation coincides with the $u$-$U$ formulation if and only if 
\begin{equation}
    d_u \subseteq \trialspacep 
    \quad \forall \, d_u \in \grad \cdot \trialspaceu
    \quad \mathrm{and} \quad
    d_U \subseteq \trialspacep
    \quad \forall \, d_U \in \grad \cdot \trialspaceU.
\end{equation}
These requirements are fulfilled when we employ these divergence-conforming spaces to discretize both displacement fields. We will demonstrate this in Sections \ref{chap:1Dcase} and \ref{chap:2Dcase}.
Note that in a multi-dimensional setting, these elements only form a suitable discretization of $\matdisp$ when $\alpha_i \geq 1$, since the trial space $\trialspaceu$ must be in $\mathbf{H}^1$. Thus, this is only possible with the higher-order continuity provided by IGA.

\begin{remark}    
    Mapping a divergence-conforming complex onto a general geometry can destroy its compatible structure. This can be avoided by mapping through the push-forward operation \cite{buffaIsogeometricDiscreteDifferential2011}. 
    This is not necessary for the remainder of this work, as we consider only linear and rectilinear domains.
\end{remark}

\section{One-dimensional case: soil column}\label{chap:1Dcase}

\begin{figure}[h!]
    \centering
    \resizebox{0.5\linewidth}{!}
    {    
        \tikzsetnextfilename{1Dcase}
\begin{tikzpicture}[x=1cm,y=1cm]

\useasboundingbox (-5, -1) rectangle (6, 9);

\def\xleftbar{-4}
\def\xrightbar{4}
\def\ytop{8}
\def\ybot{0}
\def\xleft{-1}
\def\xright{1}
\def\hatchthickness{0.25}

\definecolor{domaincol} {RGB}{204, 221, 170};

\coordinate (Ltop) at (\xleft, \ytop);
\coordinate (Lbot) at (\xleft, \ybot);
\coordinate (Rtop) at (\xright, \ytop);
\coordinate (Rbot) at (\xright, \ybot);

\coordinate (LBbar) at (\xleftbar, \ybot);
\coordinate (RBbar) at (\xrightbar, \ybot);

\coordinate (LTbar) at (\xleftbar , \ytop);
\coordinate (RTbar) at (\xrightbar, \ytop);

\fill[domaincol]
    (Ltop) -- (Rtop) -- (Rbot) -- (Lbot) -- cycle;

\draw[line width=1pt]
    (LBbar) -- (RBbar); 

\draw[line width=1pt]
    (LTbar) -- (RTbar);
\fill[pattern=north east lines] (LTbar) rectangle ($(RTbar) + (0, \hatchthickness)$) ;

\node[right] at (RTbar) {\Large$x=50$ cm};
\node[right] at (RBbar) {\Large$x=0$};

\draw[dashed] (Lbot) -- (Ltop);
\draw[dashed] (Rbot) -- (Rtop);

\draw[line width=1pt, ->]
    ( $ (Rbot) + (0.5, 0) $ ) -- ( $ (Rbot) + (0.5, 2) $ );

\coordinate (Apos) at (0.0, {(\ytop - \ybot) / 5.0 * 3.0});
\node at (Apos) {$\bullet$};
\node[anchor=south west] at ($(Apos)+(0.08,0.08)$) {\Large$A$};
\draw[-] ($(Apos)+(0.25,0)$) -- ($(Apos) + (\xrightbar, 0)$)
  node[right] {\Large$x=30$ cm};
  
\end{tikzpicture}
    }
    \caption{Schematic of the 1D benchmark. A constant velocity is instantaneously applied to one end of a soil column.}
    \label{fig:1Dcase}
\end{figure}
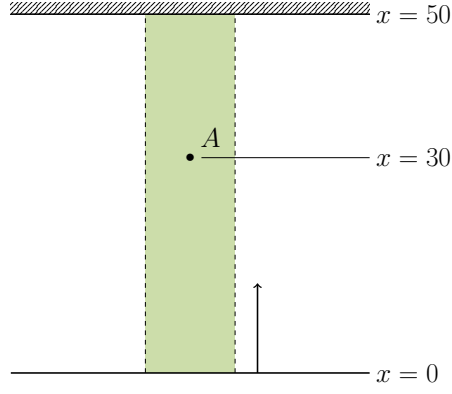

Our first case study is the 50 cm long 1D soil column previously analyzed by, among others, Gajo et al. \cite{gajoEvaluationThreeTwofield1994}. A schematic of the scenario is provided in \autoref{fig:1Dcase}. The column has unit area, and symmetry conditions along the sides. It is fixed at $x=$ \num{50} \unit{\centi\meter}, while at $x=$ \num{0} \unit{\centi\meter}, a constant positive velocity is imposed from $t=$ \num{0} \unit{\second}, resulting in a shock wave reflecting through the column. The top and bottom are impermeable. Hence, both conditions apply to both the matrix displacement $\matdisp$ and the fluid displacement $\fludisp$. The material parameters are provided in \autoref{tab:matpars_1D}. Note that the compressibility $Q$ is a resultant of these \cite{coussyPoromechanics2004}.\\
We validate our implementation using the problem defined as such in \cref{subsec:validation1D}. For our spectral analysis we consider free vibration, where the column is fixed at $x=$ \num{0} \unit{\centi\meter} and free at $x=$ \num{50} \unit{\centi\meter}. To allow normalization of spectra, we first discuss the analytical solution in \cref{subsec:analytical1D}. Since there are two branches of solutions, we discuss a method for labeling discrete modes accordingly in \cref{subsec:discreteprelims}. With that, in \crefrange{subsec:1Dspectra-uU}{subsec:1D_formulation_comparision} we study spectra for the $u$-$U$ formulation, (non)conforming spaces for the $u$-$p$-$U$ formulation, and compare the two formulations against one another. Finally, we study the critical timestep with- and without outlier removal in \cref{subsec:1D_outlierremoval,subsec:1Dadmissibletimestep}.

\begin{table}[h!]
    \centering
    \renewcommand{\arraystretch}{1.25}
    \begin{tabular}{ccc}
        \toprule
        Parameter & Value & Unit \\
        \midrule
        $\lambda$ & \num{5.15e10}  & \unit{\gram\per\centi\meter\per\second} \\
        $\mu$     & \num{9.91e10}  & \unit{\gram\per\centi\meter\per\second} \\
        $\alpha$  & \num{0.6734}   & -- \\
        $\phi$    & \num{0.18}     & -- \\
        $k$       & \num{0.148e-2} & \unit{\centi\meter\cubed\second\per\gram} \\
        $\rho_s$  & \num{2.66}     & \unit{\gram\per\centi\meter\cubed} \\
        $\rho_f$  & \num{1.0}      & \unit{\gram\per\centi\meter\cubed} \\
        \bottomrule
    \end{tabular}
    \caption{Material parameters for the 1D benchmark. The bulk modulus is $K = \lambda + 2\mu/3$, with which the skeleton bulk modulus is $K_s = \frac{K}{1-\alpha}$. Then, the inverse of the compressibility modulus is $Q^{-1} = (\alpha  - \phi) / K_s$ \cite{coussyPoromechanics2004}.}
    \label{tab:matpars_1D}
\end{table}

\subsection{Validation}\label{subsec:validation1D}

\begin{figure}[h!]
    \centering
    
    \def\marksize{2pt}
    \tikzsetnextfilename{newmark_validation}
\begin{tikzpicture}

\begin{groupplot}[
    group style={
        group size=1 by 2,
        vertical sep=1.4cm,
    },
    width=.8\textwidth,
    height=.35\textwidth,
    grid=both,
    xlabel={Time (t/T)},
    cycle list name=mytab10,
]

\nextgroupplot[
    ylabel={Solid velocity [cm/s]},
    ymin=-0.25,
    ymax=1.25,
    xmin=0,
    xmax=8,
    legend columns=2,
    legend to name=newmark_validation_legend,
    legend image post style={
        mark indices={2}
    },
    legend style={
        font=\footnotesize,
        rounded corners=2pt,
        draw=black!35,
    },
]

\addplot+[
    mark indices={0,30,...,10000},
    mark size=\marksize,
    restrict y to domain=-10:10,
    unbounded coords=discard
]
table[col sep=comma, x=time_normalized, y=solid_velocity_cm_per_s]
{data/1D_solutions_over_time/newmark_validation_ne100_p3-2-3_c2-1-2_newmark_dt1em06_x20.csv};
\addlegendentry{$u$-$p$-$U$}

\addplot+[
    mark indices={0,30,...,10000},
    mark size=\marksize,    
    restrict y to domain=-10:10,
    unbounded coords=discard
]
table[col sep=comma, x=time_normalized, y=solid_velocity_cm_per_s]
{data/1D_solutions_over_time/newmark_validation_ne100_p3-3_c2-2_newmark_dt1em06_x20.csv};
\addlegendentry{$u$-$U$}

\nextgroupplot[
    ylabel={Fluid velocity [cm/s]},
    ymin=-0.25,
    ymax=1.25,
    xmin=0,
    xmax=8,
]

\addplot+[
    mark indices={0,30,...,10000},
    mark size=\marksize,
    restrict y to domain=-10:10,
    unbounded coords=discard
]
table[col sep=comma, x=time_normalized, y=fluid_velocity_cm_per_s]
{data/1D_solutions_over_time/newmark_validation_ne100_p3-2-3_c2-1-2_newmark_dt1em06_x20.csv};

\addplot+[
    mark indices={0,30,...,10000},
    mark size=\marksize,
    restrict y to domain=-10:10,
    unbounded coords=discard
]
table[col sep=comma, x=time_normalized, y=fluid_velocity_cm_per_s]
{data/1D_solutions_over_time/newmark_validation_ne100_p3-3_c2-2_newmark_dt1em06_x20.csv};

\end{groupplot}

\end{tikzpicture}
    
    \begingroup
    \tikzexternaldisable
    \pgfplotslegendfromname{newmark_validation_legend}
    \endgroup    
    
    \vspace{0.5em}
    \caption{The velocity of the matrix material and the pore fluid at point $A$ in \autoref{fig:1Dcase} over time. Solutions for both the $u$-$p$-$U$ and the $u$-$U$ formulation. The solid- and fluid fields are discretized with cubic B-splines, the pressure field with quadratics, all at optimal regularity.}
    \label{fig:solution_over_time}
\end{figure}
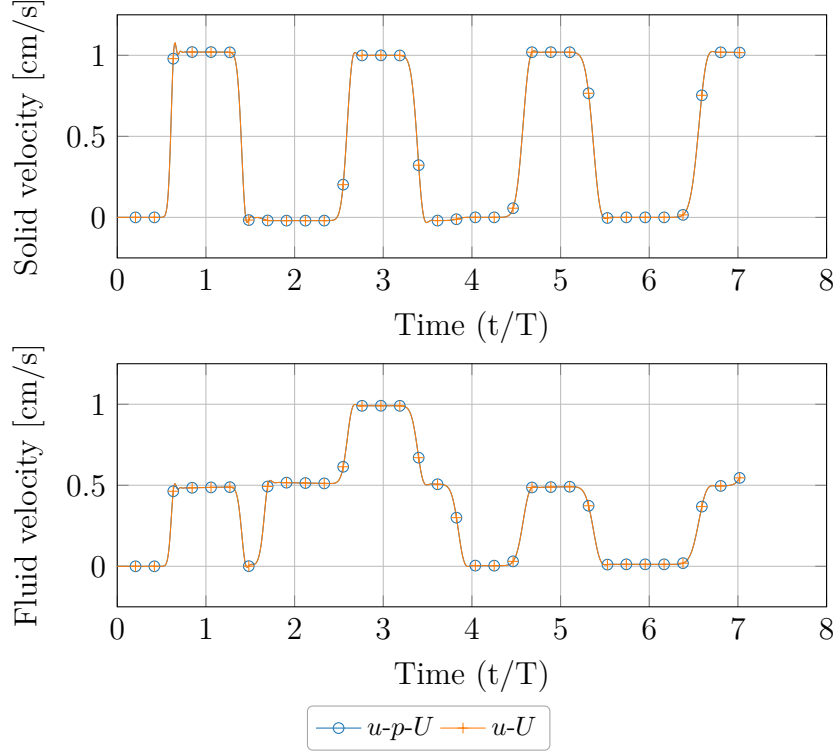

The domain is discretized with 100 evenly spaced elements (non-zero knot intervals). Following Ref. \cite{gajoEvaluationThreeTwofield1994}, the simulation runtime, \num{990e-6} \unit{\second}, is discretized using the implicit Newmark scheme. We empirically determined the Newmark parameters as $\gamma = 0.6$ and $\beta = 0.3$, and use 990 timesteps of \num{1e-6} \unit{\second}. \\
To assess the validity of the considered formulations, we solved the benchmark using both the $u$-$p$-$U$ and $u$-$U$ formulations with cubic B-splines for the displacements and quadratics for the pressure, all at optimal regularity. In \autoref{fig:solution_over_time}, the resulting velocities $\dot{\matdisp}$ and $\dot{\fludisp}$ at point $A$ in \autoref{fig:1Dcase} are plotted over time.
The time-axis is normalized by the characteristic timescale $T=$ \num{140.9} \unit{\micro\second} \cite{gajoEvaluationThreeTwofield1994}. We see one wave that is observable in both fields and a second wave which is only visible in the fluid. These correspond to the fast and slow wave discussed in \autoref{subsec:analytical1D}. The two simulations are nearly indistinguishable, and agree well with the findings reported in \cite{gajoEvaluationThreeTwofield1994}, thereby validating our implementation.

\subsection{Analytical preliminaries}\label{subsec:analytical1D}

To assess the error in our discrete eigenvalues, we compare them to analytical results. We here consider only the undamped, 1D case. A more general treatment may be found in Ref. \cite{carcioneChapter7Biot2022}. \\
Our starting point is the strong $u$-$U$ formulation \eqref{eq:twofield_strong}. Herein, we neglect damping by letting $k\rightarrow \infty$. After reducing to 1D, substituting Hooke's law \eqref{eq:linela}, and assuming homogeneous parameters, we can write the two equations as the linear system
\begin{align}
    \begin{bmatrix}
        \lambda+2\mu + \aph^2 Q
        &
        \aph Q \phi
        \\
        \aph Q \phi & Q\phi^2 
    \end{bmatrix}    
    \begin{bmatrix}
        \frac{\partial^2 u}{\partial x^2} \\ \frac{\partial^2 U}{\partial x^2}
    \end{bmatrix}
    =  
    \begin{bmatrix}
        \rho_1 & 0
        \\        
        0 & \rho_2
    \end{bmatrix}    
    \begin{bmatrix}
        \frac{\partial^2 u}{\partial t^2} \\ \frac{\partial^2 U}{\partial t^2}
    \end{bmatrix},
\end{align}
with $\rho_1 = (1-\phi)\rho_s$ and $\rho_2 = \phi \rho_f$. We then substitute the plane wave ansatz, $u = \tilde{u} \exp( i(k x - \omega t))$, and $U = \tilde{U} \exp( i(kx - \omega t))$ to obtain
\begin{align}\label{eq:analyticaleigproblem}
    \underbrace
    {    
    \left(
    k^2
    \begin{bmatrix}
        \lambda+2\mu + \aph^2 Q
        &
        \aph Q \phi
        \\
        \aph Q \phi & Q\phi^2 
    \end{bmatrix}    
    -
    \omega^2 
    \begin{bmatrix}
        \rho_1 & 0
        \\        
        0 & \rho_2
    \end{bmatrix} 
    \right)
    }_{:= \matr{D}}
    \begin{bmatrix}
        u \\ U
    \end{bmatrix}
    = \begin{bmatrix}
        0 \\ 0
    \end{bmatrix}.
\end{align}
Nontrivial solutions are only obtained if $\det(\matr{D}) = 0$. Explicitly, that is if
\begin{align}
\begin{split}
    \omega^4 \rho_1 \rho_2 - k^2 \omega^2 \underbrace{\left[(\lambda+2\mu + \aph^2 Q) \rho_2 + Q\phi^2\rho_1\right]}_{:= S} &+ \\
    k^4 (\lambda+2\mu) Q \phi^2 &= 0.
\end{split}
\end{align}
Since this is a quadratic equation in $\omega^2$, two solution branches are identified through the quadratic equation,
\begin{align}\label{eq:omega_analytical}
    \omega^2 = k^2 \underbrace{\left[\frac{1}{2\rho_1\rho_2} \left(S \pm \sqrt{S^2 -4\rho_1\rho_2(\lambda+2\mu)Q\phi^2}\right)\right]}_{:=c_\pm^2}.
\end{align}
This identifies two wave speeds, $c_+$ and $c_-$. Back-substituting $\omega^2 = k^2 c_\pm^2$ into the bottom equation in \eqref{eq:analyticaleigproblem} yields, for nonzero $k$, 
\begin{align}\label{eq:participationratios}
    \frac{u}{U} = \frac{\rho_2 c^2_\pm - Q \phi^2}{\aph Q \phi} :=r_\pm,
\end{align}
where $r_\pm$ is the mode participation ratio. \\
At this point, independent of the geometry of the problem, we have identified two kinds of waves. Namely the fast Biot wave with $c_+, r_+$ and the slow Biot wave with $c_-, r_-$. The participation ratio of the former is positive, while that of the latter is negative \cite{carcioneChapter7Biot2022}. Thus, in the fast wave the fluid oscillates in phase with the solid, while the opposite occurs in the slow wave. \\
In our discrete analysis, we will employ Dirichlet boundary conditions on one side of the domain and Neumann conditions on the other side. Starting with a harmonic ansatz, after requiring those boundary conditions to hold, we find the mode shapes
\begin{equation}{\label{eq:generalmodeshape}}
    U_n(x) = A_n\sin(k_nx), \qquad u_n(x) = r_\pm U_n,
\end{equation}
with wave numbers 
\begin{equation}
    k_n = (2n+1)\pi / (2L).
\end{equation}
Together with \eqref{eq:omega_analytical}, this defines the natural frequencies.

\subsection{Discrete preliminaries}\label{subsec:discreteprelims}
\begin{figure}
    \centering
    \def\markersize{1.2pt}

\pgfplotstableread[col sep=comma]{data/spectra_1D/spectrum_p2-2_c1-1.csv}\parttable

\pgfplotstableread[col sep=comma]{data/analytical_participation_ratios.csv}\reftable

\pgfplotstablegetelem{0}{FastRatio}\of\reftable
\edef\rfast{\pgfplotsretval}

\pgfplotstablegetelem{0}{SlowRatio}\of\reftable
\edef\rslow{\pgfplotsretval}

\tikzsetnextfilename{participation_ratio}
\begin{tikzpicture}
\begin{axis}[
  width=.8\linewidth,
  height=.5\linewidth,
  xlabel={$n$},
  ylabel={$r_n$},
  xmin=-20,
  xmax=550,
  ymin=-0.20,
  ymax=2.3,
  legend style={
    at={(0.98,0.50)},
    anchor=east,
    fill=none,
    rounded corners=2pt,
    draw=black!35,
  },
]

\addplot[
  tabblue,
  dashed,
  thick,
  mark=none
]
coordinates {(-20,\rfast) (550,\rfast)};
\addlegendentry{$r_+$}

\addplot[
  taborange,
  dashed,
  thick,
  mark=none
]
coordinates {(-20,\rslow) (550,\rslow)};
\addlegendentry{$r_-$}

\addplot[
  scatter,
  only marks,
  mark=*,
  mark size=\markersize,
  scatter src=explicit symbolic,
  scatter/classes={
    0={tabblue, mark options={solid, fill=tabblue}},
    1={taborange, mark options={solid, fill=taborange}}
  }
]
table[
  x expr=\coordindex+1,
  y=ParticipationRatio,
  meta=ModeLabel,
  col sep=comma
]{\parttable};

\addlegendimage{
  only marks,
  mark=*,
  mark size=\markersize,
  tabblue,
  mark options={solid, fill=tabblue}
}
\addlegendentry{Fast discrete modes}

\addlegendimage{
  only marks,
  mark=*,
  mark size=\markersize,
  taborange,
  mark options={solid, fill=taborange}
}
\addlegendentry{Slow discrete modes}

\end{axis}
\end{tikzpicture}
    \caption{Participation ratios of the discrete spectrum for the $u$-$U$ formulation, using 256 quadratic B-splines for both fields.}
    \label{fig:uUparticipation}
\end{figure}
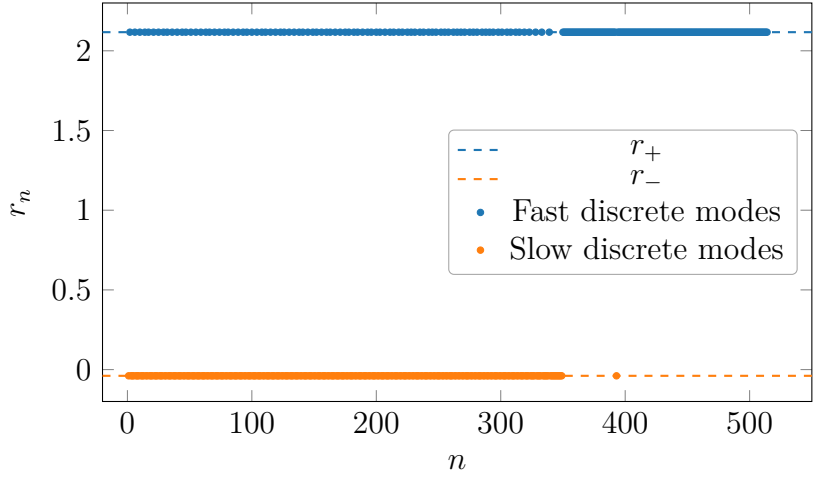

Since \eqref{eq:secondorderform} is a damped second-order system, the corresponding spectrum may be computed by solving a Quadratic Eigenvalue Problem (QEP) \cite{tisseurQuadraticEigenvalueProblem2001}. In \ref{app:damped}, we discuss this further, and solve this damped problem for one discretization. There, we observe that the current choice of parameters leads to eigenfrequencies of which the imaginary part is negligible relative to their real part. Therefore, the critical timestep is well-informed by the undamped problem. Since the damped problem introduces complexity, especially in the analytical setting, without further informing differences between discretizations, we neglect damping in the spectral analysis that follows. Consequently, the QEP simplifies to the Generalized Eigenvalue Problem (GEP)
\begin{equation}\label{eq:generalizedeigenvalueproblem}
    \overline{\matr{M}} \matr{\varphi} = \omega^2 \, \overline{\matr{K}} \matr{\varphi}, \quad \mathrm{with} \quad
    \matr{\varphi} = 
    \begin{bmatrix}
        \matr{\varphi_u}\\ \matr{\varphi_U}
    \end{bmatrix},
\end{equation}
where $\matr{\varphi}$ is the eigenvector and $\omega$ the eigenfrequency. Since this is a symmetric system, all eigenfrequencies are real \cite{golubChapter8Symmetric2013}.\\
To compare discrete eigenfrequencies to analytical eigenfrequencies, we need to categorize the discrete ones into fast and slow modes, as defined in \autoref{subsec:analytical1D}. To that end, we use the participation ratios \eqref{eq:participationratios}. We find these numerically as
\begin{equation}
    r_n = \frac
    {(u^h_n, U^h_n)_{L^2}}
    {(U^h_n, U^h_n)_{L^2}}
    ,
\end{equation}
where $(\Box, \Box)$ denotes the inner product and the subscripts $n$ the mode number. For the $u$-$U$ formulation, discretized with quadratic B-splines, these ratios are displayed in \autoref{fig:uUparticipation}. We find two clusters, each of which closely matches the analytically determined ratios. Thus, we can indeed label discrete modes as fast and slow waves according to this metric.

\subsection{Discrete spectra: the \texorpdfstring{$u$-$U$}{u-U} formulation}\label{subsec:1Dspectra-uU}
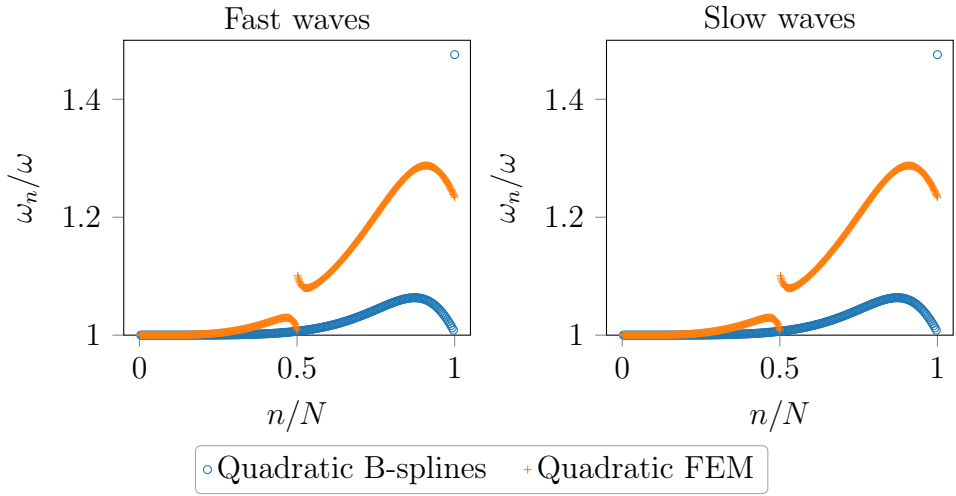
\begin{figure}[h!]
    \centering

    \def\marksize{1.5pt}

\tikzsetnextfilename{normalized_quadratic_uU}
\begin{tikzpicture}
\begin{groupplot}[
    group style={
        group size=2 by 1,
        horizontal sep=1.8cm,
    },
    width=0.45\textwidth,
    height=0.40\textwidth,
    xmin=-0.05,
    xmax=1.05,
    ymin=1.0,
    ymax=1.5,
    xlabel={$n/N$},
    ylabel={$\omega_n / \omega$},
    grid=none,
    tick align=outside,
    xtick pos=left,
    ytick pos=left,
]

\nextgroupplot[
    title={Fast waves},
    legend to name=normalized_quadratic_uU_legend,
    legend columns=2,
    legend style={
        /tikz/every even column/.append style={column sep=1em},
        rounded corners=2pt,
        draw=black!35,
    },
]

\addplot+[
    tabblueplot,
    only marks,
    mark size=\marksize,
    mark options={solid, fill=tabblue},
]
table[
    col sep=comma,
    x=NormalizedSegmentedModenumber,
    y=NormalizedFrequency,
] {data/spectra_1D_segmented/spectrum_p2-2_c1-1_fast.csv};
\addlegendentry{Quadratic B-splines}

\addplot+[
    taborangeplot,
    only marks,
    mark size=\marksize,
    mark options={solid, fill=taborange},
]
table[
    col sep=comma,
    x=NormalizedSegmentedModenumber,
    y=NormalizedFrequency,
    restrict expr to domain={\thisrow{ModeLabel}}{0:0},
] {data/spectra_1D_segmented/spectrum_p2-2_c0-0_fast.csv};
\addlegendentry{Quadratic FEM}

\nextgroupplot[
    title={Slow waves},
]

\addplot+[
    tabblueplot,
    only marks,
    mark size=\marksize,
    mark options={solid, fill=tabblue},
]
table[
    col sep=comma,
    x=NormalizedSegmentedModenumber,
    y=NormalizedFrequency,
] {data/spectra_1D_segmented/spectrum_p2-2_c1-1_slow.csv};

\addplot+[
    taborangeplot,
    only marks,
    mark size=\marksize,
    mark options={solid, fill=taborange},
]
table[
    col sep=comma,
    x=NormalizedSegmentedModenumber,
    y=NormalizedFrequency,
] {data/spectra_1D_segmented/spectrum_p2-2_c0-0_slow.csv};

\end{groupplot}
\end{tikzpicture}

    \begingroup
    \tikzexternaldisable
    \pgfplotslegendfromname{normalized_quadratic_uU_legend}
    \endgroup

    \vspace{0.5em}
    \caption{Normalized discrete spectra for the 1D poroelasticity problem, comparing quadratic FEM and B-splines. The spectra are split into slow waves and fast waves through the labeling in \autoref{fig:uUparticipation}.}
    \label{fig:1D_normalizedspectrum_quad}
\end{figure}

\begin{figure}[h!]
    \centering

    \def\marksize{1.5pt}

\tikzsetnextfilename{normalized_higherorder_uU}
\begin{tikzpicture}
\begin{groupplot}[
    group style={
        group size=2 by 1,
        horizontal sep=2.0cm, 
    },
    width=0.45\textwidth,
    height=0.40\textwidth,
    xmin=-0.05,
    xmax=1.05,
    ymin=1.0,
    ymax=1.07,
    restrict y to domain=0.9:1.1,
    xlabel={$n/N$},
    ylabel={$\omega_n / \omega$},
    grid=none,
    tick align=outside,
    xtick pos=left,
    ytick pos=left,
]

\nextgroupplot[
    title={Fast waves},
    legend to name=normalized_higherorder_uU_legend,
    legend columns=4,
    legend style={
        /tikz/every even column/.append style={column sep=1em},
        rounded corners=2pt,
        draw=black!35,
    },
]

\addplot+[
    tabblueplot,
    only marks,
    mark size=\marksize,
    mark options={solid, fill=tabblue},
]
table[
    col sep=comma,
    x=NormalizedSegmentedModenumber,
    y=NormalizedFrequency,
] {data/spectra_1D_segmented/spectrum_p2-2_c1-1_fast.csv};
\addlegendentry{Quadratic}

\addplot+[
    taborangeplot,
    only marks,
    mark size=\marksize,
    mark options={solid, fill=taborange},
]
table[
    col sep=comma,
    x=NormalizedSegmentedModenumber,
    y=NormalizedFrequency,
] {data/spectra_1D_segmented/spectrum_p3-3_c2-2_fast.csv};
\addlegendentry{Cubic}

\addplot+[
    tabgreenplot,
    only marks,
    mark size=\marksize,
    mark options={solid, fill=tabgreen},
]
table[
    col sep=comma,
    x=NormalizedSegmentedModenumber,
    y=NormalizedFrequency,
] {data/spectra_1D_segmented/spectrum_p4-4_c3-3_fast.csv};
\addlegendentry{Quartic}

\addplot+[
    tabredplot,
    only marks,
    mark size=\marksize,
    mark options={solid, fill=tabred},
]
table[
    col sep=comma,
    x=NormalizedSegmentedModenumber,
    y=NormalizedFrequency,
] {data/spectra_1D_segmented/spectrum_p5-5_c4-4_fast.csv};
\addlegendentry{Quintic}

\nextgroupplot[
    title={Slow waves},
]

\addplot+[
    tabblueplot,
    only marks,
    mark size=\marksize,
    mark options={solid, fill=tabblue},
]
table[
    col sep=comma,
    x=NormalizedSegmentedModenumber,
    y=NormalizedFrequency,
] {data/spectra_1D_segmented/spectrum_p2-2_c1-1_slow.csv};

\addplot+[
    taborangeplot,
    only marks,
    mark size=\marksize,
    mark options={solid, fill=taborange},
]
table[
    col sep=comma,
    x=NormalizedSegmentedModenumber,
    y=NormalizedFrequency,
] {data/spectra_1D_segmented/spectrum_p3-3_c2-2_slow.csv};

\addplot+[
    tabgreenplot,
    only marks,
    mark size=\marksize,
    mark options={solid, fill=tabgreen},
]
table[
    col sep=comma,
    x=NormalizedSegmentedModenumber,
    y=NormalizedFrequency,
] {data/spectra_1D_segmented/spectrum_p4-4_c3-3_slow.csv};

\addplot+[
    tabredplot,
    only marks,
    mark size=\marksize,
    mark options={solid, fill=tabred},
]
table[
    col sep=comma,
    x=NormalizedSegmentedModenumber,
    y=NormalizedFrequency,
] {data/spectra_1D_segmented/spectrum_p5-5_c4-4_slow.csv};

\end{groupplot}
\end{tikzpicture}

    \begingroup
    \tikzexternaldisable
    \pgfplotslegendfromname{normalized_higherorder_uU_legend}
    \endgroup

    \vspace{0.5em}
    \caption{Normalized discrete spectra for the poro-elasticity problem, using higher-order B-splines. Note that outliers lie outside these figures.}
    \label{fig:1D_normalizedspectrum_higherorder}
\end{figure}
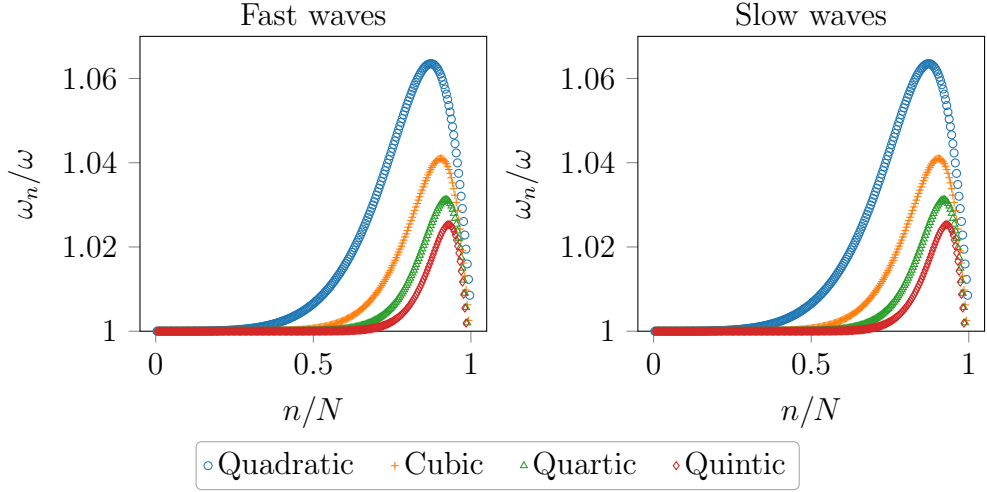

Starting with the $u$-$U$ formulation, following Ref. \cite{cottrellIsogeometricAnalysisStructural2006}, we first study the spectra resulting from quadratic Lagrangian FEM and B-splines in \autoref{fig:1D_normalizedspectrum_quad}. Note that the eigenfrequencies following from the generalized eigenvalue problem \eqref{eq:generalizedeigenvalueproblem} depend on the space, not the basis. Hence, a $C^0$ spline basis and Lagrange FEM basis of the same order result in the same eigenfrequencies. The modes in \autoref{fig:1D_normalizedspectrum_quad} are labeled, and they are normalized with respect to the analytically computed eigenfrequencies \eqref{eq:omega_analytical}.\\
We find essentially the same results as reported for linear elasticity \cite{cottrellIsogeometricAnalysisStructural2006}, but twice: once for each wave-type. That is, approximately half of the modes resulting from the FEM-discretization are non-physical optical modes, while the use of B-splines effectively eliminates these modes. Splines do, however, introduce outlier frequencies, the number of which is constant under mesh refinement. Hence, we can conclude that the entire linear system resulting from B-splines, modulo the outlier modes, contributes to finding a meaningful solution, in contrast to half of the modes in the FEM case. \\
In \autoref{fig:1D_normalizedspectrum_higherorder}, we increase the polynomial order of the B-spline basis. This too resembles twice the figure reported for elasticity in Ref. \cite{cottrellIsogeometricAnalysisStructural2006}. The entire spectrum converges toward the reference solution as the order increases, with the exception of outliers, which grow both in number and in magnitude. These results convey that the advantageous spectral properties of B-splines, as observed for linear elasticity among others, do carry over to poroelastodynamics.

\subsection{Compatible spaces}
\begin{figure}[h!]
    \centering
    \def\markersize{1.2pt}

\pgfplotstableread[col sep=comma]{data/analytical_participation_ratios.csv}\parttable
\pgfplotstablegetelem{0}{FastRatio}\of\parttable
\edef\rfast{\pgfplotsretval}
\pgfplotstablegetelem{0}{SlowRatio}\of\parttable
\edef\rslow{\pgfplotsretval}

\tikzsetnextfilename{participation_ratio_nonconforming}
\begin{tikzpicture}
\begin{axis}[
  width=.8\linewidth,
  height=.5\linewidth,
  xlabel={$n$},
  ylabel={$r_n$},
  xmin=-20,
  xmax=1100,
  ymin=-0.20,
  ymax=6.5,
  legend style={
    at={(0.03,0.71)},
    anchor=west,
    fill=none,
    rounded corners=2pt,
    draw=black!35,
  },
]

\addplot[
  tabblue,
  dashed,
  thick,
  mark=none
]
coordinates {(-20,\rfast) (1200,\rfast)};
\addlegendentry{$r_+$}

\addplot[
  taborange,
  dashed,
  thick,
  mark=none
]
coordinates {(-20,\rslow) (1200,\rslow)};
\addlegendentry{$r_-$}

\addplot[
  scatter,
  only marks,
  mark=*,
  mark size=\markersize,
  scatter src=explicit symbolic,
  scatter/classes={
    0={tabblue, mark options={solid, fill=tabblue}},
    1={taborange, mark options={solid, fill=taborange}},
    2={tabgreen, mark options={solid, fill=tabgreen}},
    3={tabgreen, mark options={solid, fill=tabgreen}}
  }
]
table[
  x expr=\coordindex+1,
  y=ParticipationRatio,
  meta=ModeLabel,
  col sep=comma
]{data/spectra_1D/spectrum_p2-1-2_c0-0-0.csv};

\addlegendimage{
  only marks,
  tabblue,
  mark=*,
  mark options={solid, fill=tabblue},
  mark size=\markersize
}
\addlegendentry{Fast discrete modes}

\addlegendimage{
  only marks,
  taborange,
  mark=*,
  mark options={solid, fill=taborange},
  mark size=\markersize
}
\addlegendentry{Slow discrete modes}

\addlegendimage{
  only marks,
  tabgreen,
  mark=*,
  mark options={solid, fill=tabgreen},
  mark size=\markersize
}
\addlegendentry{Spurious modes}

\end{axis}
\end{tikzpicture}
    \caption{Participation ratios resulting from the nonconforming $\mathbb{S}^2_0\times \mathbb{S}^1_0 \times \mathbb{S}^2_0$
    triplet for the $u$-$p$-$U$ formulation, using 256 elements.}
    \label{fig:participation_nonconforming}
\end{figure}
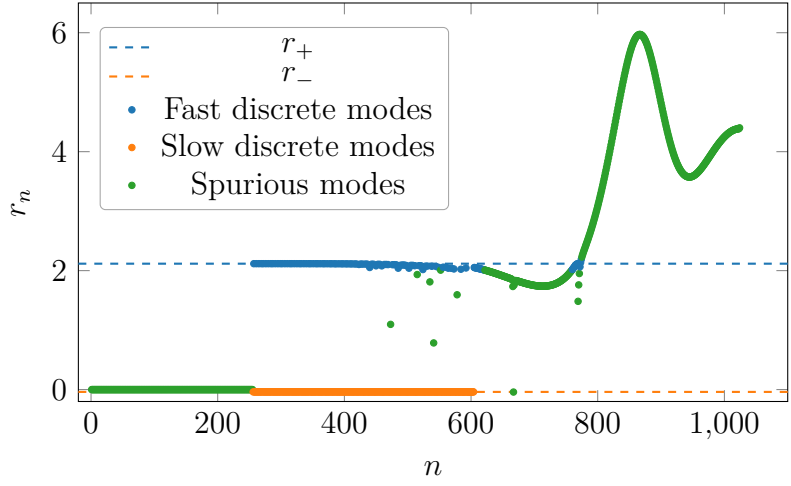

As discussed in \autoref{subsec:conforming}, a compatible discretization is required for the $u$-$p$-$U$ formulation. To illustrate the consequences of using a non-conforming discretization, we plotted the participation ratios obtained with an $\splinespace^2_0 \times \splinespace^1_0 \times \splinespace^2_0$ discretization in \autoref{fig:participation_nonconforming}. This denotes the approximation spaces used for $u$, $p$ and $U$ respectively. A similar notation is adopted for the $u$-$U$ formulation. This 1D analogue of the Taylor-Hood element clearly results in a large number of spurious eigenmodes, identifiable since they do not reproduce the analytical participation ratios. In addition, the first cluster of spurious modes has eigenvalue 0, meaning that the system does not have a unique solution. If we instead pick a conforming discretization for both $\matdisp$ and $\fludisp$, as defined in Eq. \eqref{eq:conforming1D}, we find only physical modes, reproducing the analytical participation ratios as observed in \autoref{fig:uUparticipation}.

\subsection{Comparing the \texorpdfstring{$u$-$p$-$U$}{u-p-U} and \texorpdfstring{$u$-$U$}{u-U} formulations}\label{subsec:1D_formulation_comparision}

\begin{figure}[h!]
    \centering

    \def\marksize{2.0pt}
\def\legendmarksize{4.0pt}
\def\markrepeat{24}
\def\mymarkphase{12}

\tikzsetnextfilename{normalized_twovsthree}
\begin{tikzpicture}
\begin{groupplot}[
    cycle list name=mytab10,
    group style={
        group size=2 by 1,
        horizontal sep=1.8cm,
    },
    width=0.45\textwidth,
    height=0.40\textwidth,
    xmin=-0.05,
    xmax=1.05,
    ymin=1.0,
    ymax=1.30,
    xlabel={$n/N$},
    ylabel={$\omega_n / \omega$},
    grid=none,
    tick align=outside,
    xtick pos=left,
    ytick pos=left,
    every axis plot/.append style={
        solid,
        mark repeat=\markrepeat,
        mark size=\marksize,
    },
]

\nextgroupplot[
    title={Fast waves},
    legend to name=normalized_twovsthree_legend,
    legend columns=2,
    legend style={
        /tikz/every even column/.append style={column sep=1em},
        rounded corners=2pt,
        draw=black!35,
        transpose legend,
    },
    legend image post style={
        mark size=\legendmarksize
    }
]

\addplot+ table[
    col sep=comma,
    x=NormalizedSegmentedModenumber,
    y=NormalizedFrequency,
] {data/spectra_1D_segmented/spectrum_p1-0-1_c0--1-0_fast.csv};
\addlegendentry{$\splinespace^1_0 \times \splinespace^0_{-1} \times\splinespace^1_0$}

\addplot+ table[
    col sep=comma,
    x=NormalizedSegmentedModenumber,
    y=NormalizedFrequency,
] {data/spectra_1D_segmented/spectrum_p1-1_c0-0_fast.csv};
\addlegendentry{$\splinespace^1_0 \times \splinespace^1_0$}

\addplot+[
    mark phase=\mymarkphase,
] table[
    col sep=comma,
    x=NormalizedSegmentedModenumber,
    y=NormalizedFrequency,
] {data/spectra_1D_segmented/spectrum_p2-1-2_c0--1-0_fast.csv};
\addlegendentry{$\splinespace^2_0 \times \splinespace^1_{-1} \times \splinespace^2_0$}

\addplot+[
    mark phase=\mymarkphase,
] table[
    col sep=comma,
    x=NormalizedSegmentedModenumber,
    y=NormalizedFrequency,
] {data/spectra_1D_segmented/spectrum_p2-2_c0-0_fast.csv};
\addlegendentry{$\splinespace^2_0 \times \splinespace^2_0$}

\addplot+ table[
    col sep=comma,
    x=NormalizedSegmentedModenumber,
    y=NormalizedFrequency,
] {data/spectra_1D_segmented/spectrum_p2-1-2_c1-0-1_fast.csv};
\addlegendentry{$\splinespace^2_1 \times \splinespace^1_0\times \splinespace^2_1$}

\addplot+ table[
    col sep=comma,
    x=NormalizedSegmentedModenumber,
    y=NormalizedFrequency,
] {data/spectra_1D_segmented/spectrum_p2-2_c1-1_fast.csv};
\addlegendentry{$\splinespace^2_1 \times \splinespace^2_1$}

\addplot+[
    mark phase=\mymarkphase,
] table[
    col sep=comma,
    x=NormalizedSegmentedModenumber,
    y=NormalizedFrequency,
] {data/spectra_1D_segmented/spectrum_p2-1-2_c0-0-1_fast.csv};
\addlegendentry{$\splinespace^2_0 \times \splinespace^1_0 \times \splinespace^2_1$}

\addplot+[
    mark phase=\mymarkphase,
] table[
    col sep=comma,
    x=NormalizedSegmentedModenumber,
    y=NormalizedFrequency,
] {data/spectra_1D_segmented/spectrum_p2-2_c0-1_fast.csv};
\addlegendentry{$\splinespace^2_0 \times \splinespace^2_1$}

\nextgroupplot[
    title={Slow waves},
]

\addplot+ table[
    col sep=comma,
    x=NormalizedSegmentedModenumber,
    y=NormalizedFrequency,
] {data/spectra_1D_segmented/spectrum_p1-0-1_c0--1-0_slow.csv};

\addplot+ table[
    col sep=comma,
    x=NormalizedSegmentedModenumber,
    y=NormalizedFrequency,
] {data/spectra_1D_segmented/spectrum_p1-1_c0-0_slow.csv};

\addplot+[
    mark phase=\mymarkphase,
] table[
    col sep=comma,
    x=NormalizedSegmentedModenumber,
    y=NormalizedFrequency,
] {data/spectra_1D_segmented/spectrum_p2-1-2_c0--1-0_slow.csv};

\addplot+[
    mark phase=\mymarkphase,
] table[
    col sep=comma,
    x=NormalizedSegmentedModenumber,
    y=NormalizedFrequency,
] {data/spectra_1D_segmented/spectrum_p2-2_c0-0_slow.csv};

\addplot+ table[
    col sep=comma,
    x=NormalizedSegmentedModenumber,
    y=NormalizedFrequency,
] {data/spectra_1D_segmented/spectrum_p2-1-2_c1-0-1_slow.csv};

\addplot+ table[
    col sep=comma,
    x=NormalizedSegmentedModenumber,
    y=NormalizedFrequency,
] {data/spectra_1D_segmented/spectrum_p2-2_c1-1_slow.csv};

\addplot+[
    mark phase=\mymarkphase,
] table[
    col sep=comma,
    x=NormalizedSegmentedModenumber,
    y=NormalizedFrequency,
] {data/spectra_1D_segmented/spectrum_p2-1-2_c0-0-1_slow.csv};

\addplot+[
    mark phase=\mymarkphase,
] table[
    col sep=comma,
    x=NormalizedSegmentedModenumber,
    y=NormalizedFrequency,
] {data/spectra_1D_segmented/spectrum_p2-2_c0-1_slow.csv};

\end{groupplot}
\end{tikzpicture}

    \begingroup
    \tikzexternaldisable
    \pgfplotslegendfromname{normalized_twovsthree_legend}
    \endgroup

    \vspace{0.5em}
    \caption{Discrete spectra for the $u$-$p$-$U$ and $u$-$U$ formulations using various spaces.}
    \label{fig:twovsthree}
\end{figure}
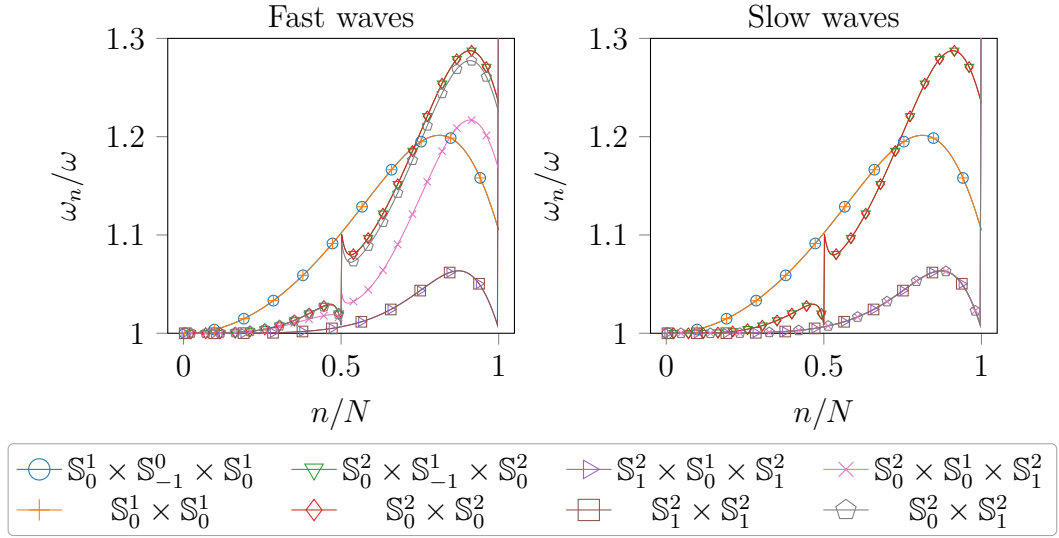

In \autoref{fig:twovsthree}, we study how the $u$-$p$-$U$ and $u$-$U$ formulations compare in terms of spectral behavior. We compare them using compatible bases which are linear, optimal regularity quadratic, and suboptimal regularity quadratic. In addition, we consider a quadratic basis wherein only the basis for the solid displacement $\matr{N}^u$ has suboptimal regularity. Note that all combinations satisfy the compatibility requirements \eqref{eq:conforming1D}. \\
These results confirm the expectations from \autoref{subsec:conforming}: with divergence-conforming spaces, the two formulations are equivalent. Hence, they result in exactly the same spectrum. Conversely, while the discretizations where the solid displacement is not divergence-conforming appear stable, they are not equivalent. In particular, the optical branches do not coincide. Note that in this last case, the optical branches do not reproduce the analytical participation ratios. Hence, classifying them as fast waves, since their $r_n$ is positive, is not as strong as for the other cases.

\subsection{The maximum eigenfrequency and outlier removal}\label{subsec:1D_outlierremoval}

\begin{figure}[h!]
    \centering

    \begin{subfigure}[t]{0.48\textwidth}
        \centering
            \def\marksize{2.0pt}
\def\markrepeat{32}

\tikzsetnextfilename{nonnormalized_spectra_outliers}
\begin{tikzpicture}
\begin{axis}[
    cycle list name=mytab10,
    width=.70\textwidth,
    height=.70\textwidth,
    scale only axis,
    xlabel={$n$ [-]},
    ylabel={$\omega_n$ [Mrad/s]},
    grid=none,
    tick align=outside,
    xtick pos=left,
    ytick pos=left,
    legend columns=3,
    legend style={
        at={(0.5,-0.30)},
        anchor=north,
        /tikz/every even column/.append style={column sep=1em},
        rounded corners=2pt,
        draw=black!35,
        transpose legend,
    },
    every axis plot/.append style={
        solid,
        mark repeat=\markrepeat,
        mark size=\marksize,
    },
]

\addplot+ table[
    col sep=comma,
    x expr=\coordindex+1,
    y expr=\thisrow{Frequency}/1e6,
] {data/spectra_1D/spectrum_p1-1_c0-0.csv};
\addlegendentry{$\splinespace^1_0 \times \splinespace^1_0$}

\pgfmathtruncatemacro{\phaseA}{mod(514,\markrepeat)}
\addplot+[
    mark phase = \phaseA,
] table[
    col sep=comma,
    x expr=\coordindex+1,
    y expr=\thisrow{Frequency}/1e6,
] {data/spectra_1D/spectrum_p2-2_c1-1.csv};
\addlegendentry{$\splinespace^2_1 \times \splinespace^2_1$}

\pgfmathtruncatemacro{\phaseB}{mod(516,\markrepeat)}
\addplot+[
    mark phase = \phaseB,
] table[
    col sep=comma,
    x expr=\coordindex+1,
    y expr=\thisrow{Frequency}/1e6,
] {data/spectra_1D/spectrum_p3-3_c2-2.csv};
\addlegendentry{$\splinespace^3_2 \times \splinespace^3_2$}

\addplot+ table[
    col sep=comma,
    x expr=\coordindex+1,
    y expr=\thisrow{Frequency}/1e6,
] {data/spectra_1D/spectrum_p2-2_c0-0.csv};
\addlegendentry{$\splinespace^2_0 \times \splinespace^2_0$}

\addplot+ table[
    col sep=comma,
    x expr=\coordindex+1,
    y expr=\thisrow{Frequency}/1e6,
] {data/spectra_1D/spectrum_p3-3_c0-0.csv};
\addlegendentry{$\splinespace^3_0 \times \splinespace^3_0$}

\draw[
    black,
    thick,
] (axis cs:-50,-0.75) rectangle (axis cs:550,7);

\end{axis}
\end{tikzpicture}
            \vfill
            \caption{Various spectra, including outliers}
            \label{fig:nonnormalized_withoutliers}
        \end{subfigure}
        \hfill
        \begin{subfigure}[t]{0.48\textwidth}
        \centering
            \def\marksize{2.0pt}
\def\markrepeat{24}

\tikzsetnextfilename{nonnormalized_spectra_outlierfree}
\begin{tikzpicture}
\begin{axis}[
    cycle list name=mytab10,
    width=.70\textwidth,
    height=.70\textwidth,
    scale only axis,
    xlabel={$n$ [-]},
    ylabel={$\omega_n$ [Mrad/s]},
    grid=none,
    tick align=outside,
    xtick pos=left,
    ytick pos=left,
    legend columns=1,
    legend style={
        at={(0.5,-0.30)},
        anchor=north,
        rounded corners=2pt,
        draw=black!35,
    },
    every axis plot/.append style={
        solid,
        mark repeat=\markrepeat,
        mark size=\marksize,
    },
    xmin=-50,
    xmax=550,
    ymin = -0.75,
    ymax=7,
]

\addplot+ table[
    col sep=comma,   
    x expr=\coordindex+1,
    y expr=\thisrow{Frequency}/1e6,
] {data/spectra_1D/spectrum_p1-1_c0-0_outlierfree.csv};
\addlegendentry{$\splinespace^1_0 \times \splinespace^1_0$}

\addplot+ table[
    col sep=comma,
    x expr=\coordindex+1,
    y expr=\thisrow{Frequency}/1e6,
] {data/spectra_1D/spectrum_p2-2_c1-1_outlierfree.csv};
\addlegendentry{$\splinespace^2_1 \times \splinespace^2_1$}

\addplot+ table[
    col sep=comma,
    x expr=\coordindex+1,
    y expr=\thisrow{Frequency}/1e6,
] {data/spectra_1D/spectrum_p3-3_c2-2_outlierfree.csv};
\addlegendentry{$\splinespace^3_2 \times \splinespace^3_2$}

\end{axis}
\end{tikzpicture}
            \vfill
            \caption{Optimal regularity only, outliers eliminated}
            \label{fig:nonnormalized_withoutoutliers}
    \end{subfigure}

    \caption{Non-normalized spectra resulting from cubic, quadratic and linear bases, with 256 elements. The two clusters correspond to fast and slow modes. In~(a) the unmodified spectra are shown. In~(b) spectra are shown from which the outliers were eliminated. For visual clarity, the spectra corresponding to $C^0$ bases are not included here. The size of this graph corresponds to the box in (a).}
    \label{fig:nonnormalized}
\end{figure}
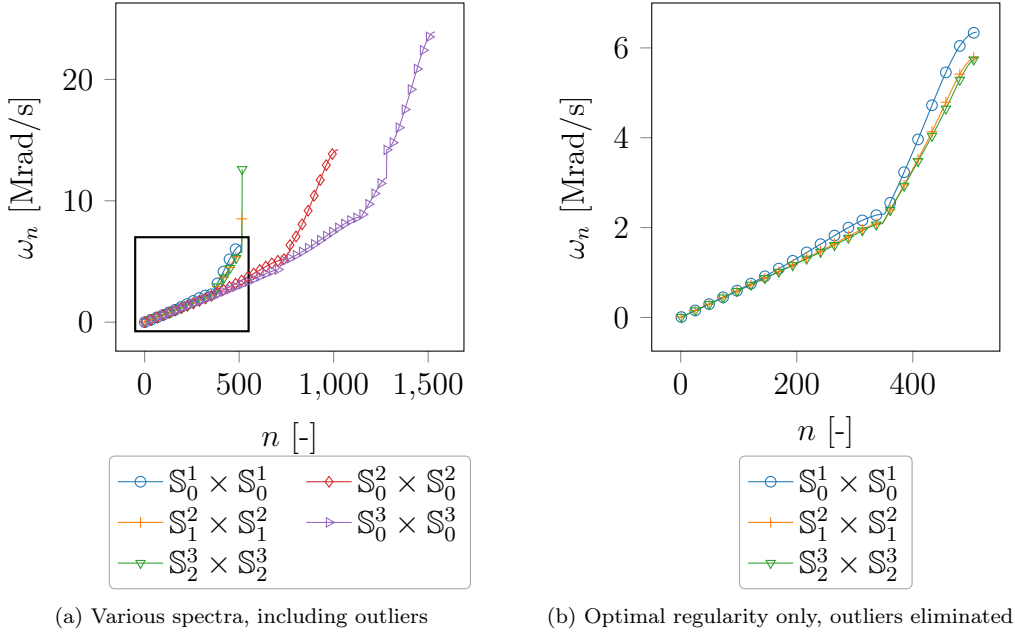

To further investigate the consequences of the choice of discretization on the critical timestep, we study spectra that are not normalized in \autoref{fig:nonnormalized_withoutliers}, so that the maximum eigenfrequency becomes visible. We observe that in both IGA and FEM the maximum frequency scales with the polynomial order. In the FEM-case the main scaling mechanism is mode count. An order elevation results in an additional number of modes equal to the number of elements. Thus, the maximum eigenfrequency scales approximately linearly with the order. Meanwhile, using optimal-regularity splines, the maximum frequency scales only with the order through the outlier modes, which grow substantially upon order elevation.\\
This unfavorable scaling motivates outlier elimination, as outlined in \autoref{sec:outliers}. A very limited number of modes can be eliminated to significantly reduce the maximum discrete frequency. We applied outlier removal to the optimal-regularity spectra in \autoref{fig:nonnormalized_withoutliers}, and present the resulting spectra in \autoref{fig:nonnormalized_withoutoutliers}. We now observe that the entire spectrum converges upon order elevation, without producing new spurious modes. The dependence of the maximum eigenfrequency on the order is now minimal, and even inverted: an order-increase decreases the maximum frequency as we converge toward a physical eigenfrequency from above. These results convey that outlier elimination prevents unfavorable scaling of the maximum eigenfrequency upon order elevation. This is a fundamental advantage of optimal regularity, outlier-free spline spaces over Lagrange FEM, which can be leveraged in the context of time integration.

\subsection{The admissible timestep}\label{subsec:1Dadmissibletimestep}

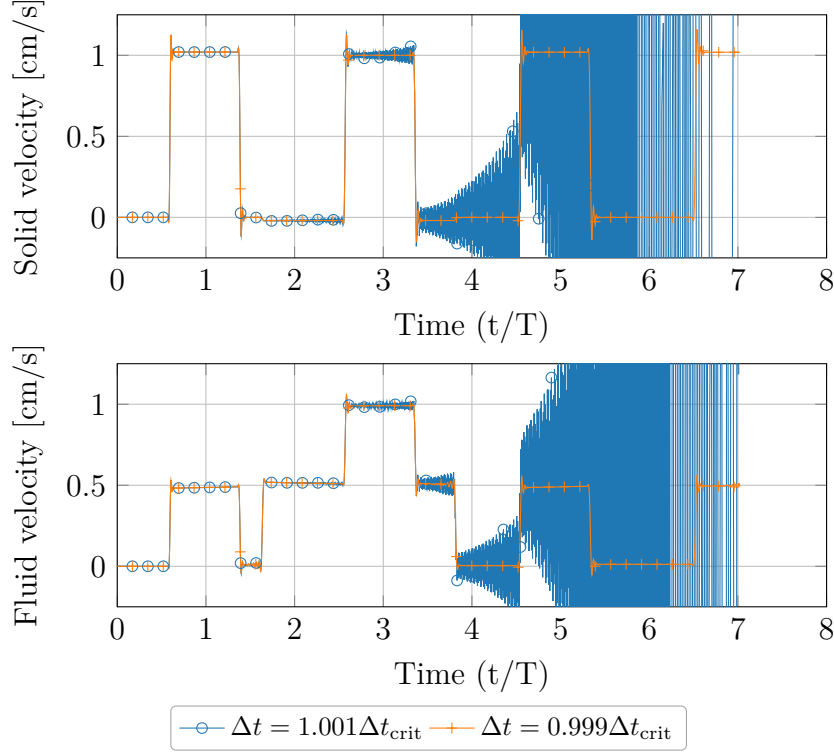
\begin{figure}[h!]
    \centering
    
    \def\marksize{2pt}
    \tikzsetnextfilename{stability_check}
\begin{tikzpicture}

\begin{groupplot}[
    group style={
        group size=1 by 2,
        vertical sep=1.4cm,
    },
    width=.8\textwidth,
    height=.35\textwidth,
    grid=both,
    xlabel={Time (t/T)},
    cycle list name=mytab10,
]

\nextgroupplot[
    ylabel={Solid velocity [cm/s]},
    ymin=-0.25,
    ymax=1.25,
    xmin=0,
    xmax=8,
    legend columns=2,
    legend to name=stability_check_legend,    
    legend image post style={
        mark indices={2}
    },    
    legend style={
        font=\footnotesize,
        rounded corners=2pt,
        draw=black!35,
    },
]

\addplot+[
    mark indices={0,50,...,10000},
    mark size=\marksize,
    restrict y to domain=-10:10,
    unbounded coords=discard
]
table[col sep=comma, x=time_normalized, y=solid_velocity_cm_per_s]
{data/1D_solutions_over_time/stability_check_ne256_p3-3_c2-2_outlierfree_rk4_dt4p913699em07_x20.csv};
\addlegendentry{$\Delta t = 1.001\Delta t_\text{crit}$}

\addplot+[
    mark indices={0,50,...,10000},
    mark size=\marksize,    
    restrict y to domain=-10:10,
    unbounded coords=discard
]
table[col sep=comma, x=time_normalized, y=solid_velocity_cm_per_s]
{data/1D_solutions_over_time/stability_check_ne256_p3-3_c2-2_outlierfree_rk4_dt4p903881em07_x20.csv};
\addlegendentry{$\Delta t = 0.999\Delta t_\text{crit}$}

\nextgroupplot[
    ylabel={Fluid velocity [cm/s]},
    ymin=-0.25,
    ymax=1.25,
    xmin=0,
    xmax=8,
]

\addplot+[
    mark indices={0,50,...,10000},
    mark size=\marksize,
    restrict y to domain=-10:10,
    unbounded coords=discard
]
table[col sep=comma, x=time_normalized, y=fluid_velocity_cm_per_s]
{data/1D_solutions_over_time/stability_check_ne256_p3-3_c2-2_outlierfree_rk4_dt4p913699em07_x20.csv};

\addplot+[
    mark indices={0,50,...,10000},
    mark size=\marksize,
    restrict y to domain=-10:10,
    unbounded coords=discard
]
table[col sep=comma, x=time_normalized, y=fluid_velocity_cm_per_s]
{data/1D_solutions_over_time/stability_check_ne256_p3-3_c2-2_outlierfree_rk4_dt4p903881em07_x20.csv};

\end{groupplot}

\end{tikzpicture}
    
    \begingroup
    \tikzexternaldisable
    \pgfplotslegendfromname{stability_check_legend}
    \endgroup    
    
    \vspace{0.5em}
    \caption{Two solutions with the same spatial discretization. One with a timestep slightly above the critical timestep, the other with a timestep slightly below. Both velocities are reported at point $A$ in \autoref{fig:1Dcase}.}
    \label{fig:stability_check}
\end{figure}

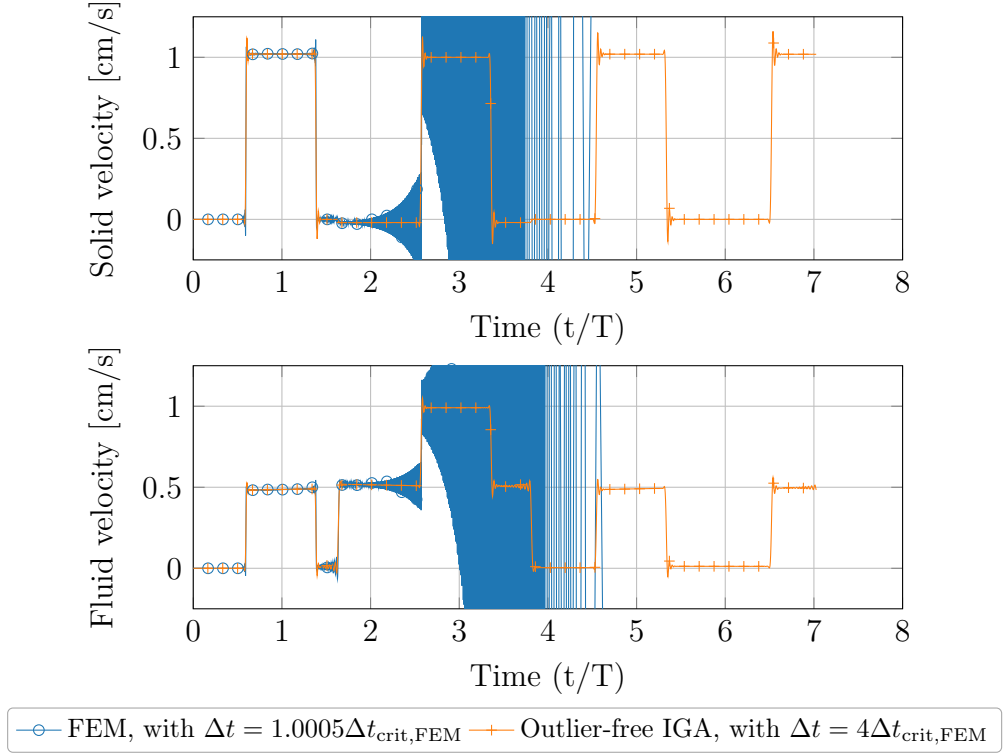
\begin{figure}[h!]
    \centering
    
    \def\marksize{2pt}
    \tikzsetnextfilename{IGA_FEM_critstep}
\begin{tikzpicture}

\begin{groupplot}[
    group style={
        group size=1 by 2,
        vertical sep=1.4cm,
    },
    width=.8\textwidth,
    height=.35\textwidth,
    grid=both,
    xlabel={Time (t/T)},
    cycle list name=mytab10,
]

\nextgroupplot[
    ylabel={Solid velocity [cm/s]},
    ymin=-0.25,
    ymax=1.25,
    xmin=0,
    xmax=8,
    legend columns=2,
    legend to name=critstepfemigalegend,    
    legend image post style={
        mark indices={2}
    },
    legend style={
        font=\footnotesize,
        rounded corners=2pt,
        draw=black!35,
    },
]

\addplot+[
    mark indices={0,200,...,10000},
    mark size=\marksize,
    restrict y to domain=-10:10,
    unbounded coords=discard
]
table[col sep=comma, x=time_normalized, y=solid_velocity_cm_per_s]
{data/1D_solutions_over_time/stability_fem_vs_iga_ne256_p3-3_c0-0_rk4_dt1p183576em07_x20.csv};
\addlegendentry{FEM, with $\Delta t = 1.0005 \Delta t_\text{crit,FEM}$}

\addplot+[
    mark indices={0,50,...,10000},
    mark size=\marksize,    
    restrict y to domain=-10:10,
    unbounded coords=discard
]
table[col sep=comma, x=time_normalized, y=solid_velocity_cm_per_s]
{data/1D_solutions_over_time/stability_fem_vs_iga_ne256_p3-3_c2-2_outlierfree_rk4_dt4p73194em07_x20.csv};
\addlegendentry{Outlier-free IGA, with $\Delta t = 4 \Delta t_\text{crit,FEM}$}

\nextgroupplot[
    ylabel={Fluid velocity [cm/s]},
    ymin=-0.25,
    ymax=1.25,
    xmin=0,
    xmax=8,
]

\addplot+[
    mark indices={0,200,...,10000},
    mark size=\marksize,
    restrict y to domain=-10:10,
    unbounded coords=discard
]
table[col sep=comma, x=time_normalized, y=fluid_velocity_cm_per_s]
{data/1D_solutions_over_time/stability_fem_vs_iga_ne256_p3-3_c0-0_rk4_dt1p183576em07_x20.csv};

\addplot+[
    mark indices={0,50,...,10000},
    mark size=\marksize,
    restrict y to domain=-10:10,
    unbounded coords=discard
]
table[col sep=comma, x=time_normalized, y=fluid_velocity_cm_per_s]
{data/1D_solutions_over_time/stability_fem_vs_iga_ne256_p3-3_c2-2_outlierfree_rk4_dt4p73194em07_x20.csv};

\end{groupplot}

\end{tikzpicture}
    
    \begingroup
    \tikzexternaldisable
    \pgfplotslegendfromname{critstepfemigalegend}
    \endgroup    
    
    \vspace{0.5em}
    \caption{Two solutions over time. One solution obtained using FEM, with a timestep just over the critical timestep. The other solution was obtained with IGA, using a timestep four times the critical timestep for FEM. Both velocities are reported at point $A$ in \autoref{fig:1Dcase}.}
    \label{fig:stability_femvsiga}
\end{figure}

\begin{figure}[h!]
    \centering

    \def\marksize{2pt}
    \tikzsetnextfilename{critstepscaling}
\begin{tikzpicture}

\begin{groupplot}[
    group style={
        group size=3 by 1,
        horizontal sep=0.3cm,
    },
    width=.38\textwidth,
    height=.38\textwidth,
    grid=none,
    xmode=log,
    ymode=log,
    xmin=10,
    xmax=700,
    ymin=0.1,
    ymax=2.0,
    xlabel={Number of elements},
    cycle list name=mytab10,
    tick align=outside,
    tick pos=left,
]

\nextgroupplot[
    title={FEM},
    ylabel={$\Delta t_\mathrm{crit} / \Delta t_\mathrm{crit,linear}$},
    legend columns=5,
    legend to name=critstepscalinglegend,
    legend image post style={
        mark indices={2}
    },
    legend style={
        rounded corners=2pt,
        draw=black!35,
    },
]

\addplot[
    dashed,
    black,
    forget plot,
] coordinates {(6,1.0) (700,1.0)};

\addplot+[
    mark size=\marksize,
]
table[col sep=comma, x=element_count, y=quadratic_normalized_by_linear]
{data/maxeigfreqinverse/maxeigfreqinverse_fem.csv};
\addlegendentry{Quadratic}

\addplot+[
    mark size=\marksize,
]
table[col sep=comma, x=element_count, y=cubic_normalized_by_linear]
{data/maxeigfreqinverse/maxeigfreqinverse_fem.csv};
\addlegendentry{Cubic}

\addplot+[
    mark size=\marksize,
]
table[col sep=comma, x=element_count, y=quartic_normalized_by_linear]
{data/maxeigfreqinverse/maxeigfreqinverse_fem.csv};
\addlegendentry{Quartic}

\addplot+[
    mark size=\marksize,
]
table[col sep=comma, x=element_count, y=quintic_normalized_by_linear]
{data/maxeigfreqinverse/maxeigfreqinverse_fem.csv};
\addlegendentry{Quintic}

\nextgroupplot[
    title={IGA},
    yticklabels={},
]

\addplot[
    dashed,
    black,
    forget plot,
] coordinates {(6,1.0) (700,1.0)};

\addplot+[
    mark size=\marksize,
]
table[col sep=comma, x=element_count, y=quadratic_normalized_by_linear]
{data/maxeigfreqinverse/maxeigfreqinverse_iga.csv};

\addplot+[
    mark size=\marksize,
]
table[col sep=comma, x=element_count, y=cubic_normalized_by_linear]
{data/maxeigfreqinverse/maxeigfreqinverse_iga.csv};

\addplot+[
    mark size=\marksize,
]
table[col sep=comma, x=element_count, y=quartic_normalized_by_linear]
{data/maxeigfreqinverse/maxeigfreqinverse_iga.csv};

\addplot+[
    mark size=\marksize,
]
table[col sep=comma, x=element_count, y=quintic_normalized_by_linear]
{data/maxeigfreqinverse/maxeigfreqinverse_iga.csv};

\nextgroupplot[
    title={outlier-free IGA},
    yticklabels={},
]

\addplot[
    dashed,
    black,
    forget plot,
] coordinates {(6,1.0) (700,1.0)};

\addplot+[
    mark size=\marksize,
]
table[col sep=comma, x=element_count, y=quadratic_normalized_by_linear]
{data/maxeigfreqinverse/maxeigfreqinverse_outlierfree_iga.csv};

\addplot+[
    mark size=\marksize,
]
table[col sep=comma, x=element_count, y=cubic_normalized_by_linear]
{data/maxeigfreqinverse/maxeigfreqinverse_outlierfree_iga.csv};

\addplot+[
    mark size=\marksize,
]
table[col sep=comma, x=element_count, y=quartic_normalized_by_linear]
{data/maxeigfreqinverse/maxeigfreqinverse_outlierfree_iga.csv};

\addplot+[
    mark size=\marksize,
]
table[col sep=comma, x=element_count, y=quintic_normalized_by_linear]
{data/maxeigfreqinverse/maxeigfreqinverse_outlierfree_iga.csv};

\end{groupplot}

\end{tikzpicture}
    \vspace{0.05em}

    \begingroup
    \tikzexternaldisable
    \pgfplotslegendfromname{critstepscalinglegend}
    \endgroup

    \vspace{0.5em}
    \caption{The critical timestep relative to that for linear FEM, for the $u$-$U$ formulation discretized with varying numbers and types of element.}    
    \label{fig:critstepscaling}
\end{figure}
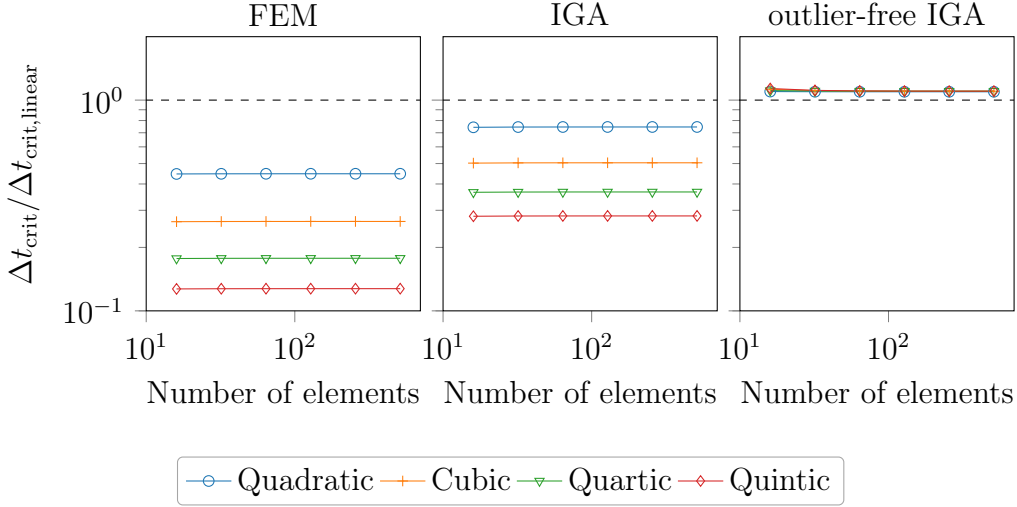

To conclude our 1D study, we illustrate the effect of discretization choice on the critical timestep in an explicit solver. The critical timestep in an explicit solver is computed through
\begin{equation}
    \Delta t_\mathrm{crit} = \frac{C_c}{\omega^h_{\mathrm{max}}},
\end{equation}
where the constant $C_c$ is the Courant number, which depends on the method. For the fourth-order Runge-Kutta method (RK4), $C_c = $ \num{2.8284} \cite{butcherNumericalMethodsOrdinary2016}.\\
We discretize our problem using the $u$-$U$ formulation with 256 elements, cubic optimal-regularity splines and outlier-removal. This results in a maximum eigenfrequency of $\omega^h_\mathrm{max} =$ \num{5.76e6} \unit{\radian}/\unit{\second}, and therefore, for RK4, a critical timestep of $\Delta t_\mathrm{crit} = $ \num{4.91e-7} \unit{\second}. In \autoref{fig:stability_check} we integrate through time using a timestep just above and just below this critical timestep.\\
As expected, a timestep just below the critical one results in a stable simulation, while a timestep just above builds up error over time at an exponential rate, confirming the validity of the spectral analyses described above. Note also that this confirms one of the observations in \ref{app:damped}: in this setting, we can compute the critical timestep in the damped problem by solving the undamped eigenvalue problem, as the largest eigenvalue has a negligible imaginary part. \\
If we use cubic FEM instead, while keeping the element count constant, the maximum eigenfrequency is \num{2.39e7} \unit{\radian}/\unit{\second}, resulting in a critical timestep for RK4 of $\Delta t_\mathrm{crit, FEM}=$ \num{1.18e-7} \unit{\second}. This is illustrated in \autoref{fig:stability_femvsiga}. Here, one solution is shown that is computed with a timestep just above this critical timestep.  As expected, this is unstable. Meanwhile, an outlier-free IGA discretization using a four times larger timestep remains stable. This illustrates that outlier-free IGA admits significantly larger timesteps than FEM.\\
In \autoref{fig:critstepscaling}, we show how the critical timestep scales with various discretization strategies. With FEM, order-elevations strongly decrease the admissible timestep, where quintics require a time step nearly ten times smaller than linears. IGA without outlier removal already scales more favorably. When outliers are removed, this results in a critical timestep that is virtually independent of the polynomial order, and even increases slightly as we converge to the maximum eigenvalue from above. This mirrors the findings for elasticity in Ref. \cite{hiemstraRemovalSpuriousOutlier2021}.

\section{Two-dimensional case: Mandel's problem}\label{chap:2Dcase}

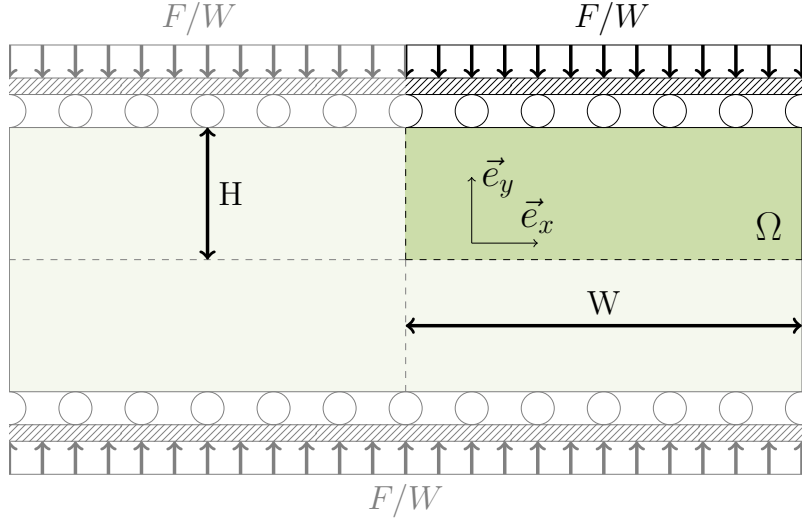
\begin{figure}[h!]
    \centering
    \resizebox{0.8\linewidth}{!}
    {    
        \tikzsetnextfilename{2Dcase}
\begin{tikzpicture}[x=1cm,y=1cm]

\useasboundingbox (-6, -4) rectangle (6, 4);

\def\xleft{-6}
\def\xright{6}
\def\ybot{-2}
\def\ytop{2}

\definecolor{domaincol} {RGB}{204, 221, 170};

\coordinate (Ltop) at (\xleft, \ytop);
\coordinate (Lbot) at (\xleft, \ybot);
\coordinate (Rtop) at (\xright, \ytop);
\coordinate (Rbot) at (\xright, \ybot);

\coordinate (Rmid) at (\xright, 0);
\coordinate (Tmid) at (0, \ytop);
\coordinate (Lmid) at (\xleft, 0);
\coordinate (Bmid) at (0, \ybot);

\coordinate (mid) at (0, 0);

\filldraw[fill=domaincol!20, draw=none] (Lbot) rectangle (Rtop);
\draw[gray] (Tmid) -- (Ltop);
\draw[gray] (Lbot) -- (Ltop);
\draw[gray] (Lbot) -- (Rbot);
\draw[gray] (Rbot) -- (Rmid);

\draw[gray, dashed] (mid) -- (Lmid);
\draw[gray, dashed] (mid) -- (Bmid);

\draw[fill=domaincol, draw=none] (mid) rectangle (Rtop);
\draw[dashed] (mid) -- (Tmid);
\draw[dashed] (mid) -- (Rmid);
\draw (Tmid) -- (Rtop);
\draw (Rtop) -- (Rmid);

\def\circrad{0.25}
\def\circledist{1}
\def\hatchedthickness{.25}

\coordinate (hatchBmid) at (0, \ytop+2*\circrad);
\coordinate (hatchTmid) at (0, \ytop+2*\circrad+\hatchedthickness);

\coordinate (BhatchTR) at (\xright+\circrad, \ybot-2*\circrad);
\coordinate (BhatchBL) at (\xleft -\circrad, \ybot-2*\circrad-\hatchedthickness);
\coordinate (BhatchBmid) at (0, \ybot-2*\circrad-\hatchedthickness);

\fill[pattern=north east lines] (hatchBmid) rectangle ($ (hatchBmid) + (\xright + \circrad, \hatchedthickness) $);
\fill[pattern=north east lines, pattern color=gray] (0, \ytop+2*\circrad) rectangle (\xleft-\circrad, \ytop+2*\circrad+\hatchedthickness);

\draw       (hatchBmid) -- ++(\xright + \circrad, 0);
\draw       (hatchTmid) -- ++(\xright + \circrad, 0);
\draw       ($ (hatchBmid) + (\xright + \circrad, 0) $) -- ++(0, \hatchedthickness);
\draw[gray] (hatchBmid) -- ++(\xleft  - \circrad, 0);
\draw[gray] (hatchTmid) -- ++(\xleft  - \circrad, 0);
\draw[gray] ($ (hatchBmid) + (\xleft - \circrad, 0) $) -- ++(0, \hatchedthickness);

\fill[pattern=north east lines, pattern color=gray] (BhatchTR) rectangle (BhatchBL);
\draw[gray] (BhatchTR) -- 
          ++(-\xright +\xleft - 2 * \circrad, 0) -- 
          ++(0, -\hatchedthickness) --
          ++(\xright -\xleft + 2 * \circrad, 0) -- 
          cycle;


\begin{scope}
  \clip (0, \ybot) rectangle (\xright, \ytop+2*\circrad);
  \draw (0, \ytop + \circrad) circle (\circrad);
\end{scope}
\begin{scope}
  \clip (\xleft, \ybot) rectangle (0, \ytop+2*\circrad);
  \draw[gray] (0, \ytop + \circrad) circle (\circrad);
\end{scope}

\draw[gray] (0, \ybot - \circrad) circle (\circrad);

\foreach \i in {1,...,6} {
    \draw       ( \i*4 * \circrad, \ytop+\circrad) circle (\circrad);
    \draw[gray] (-\i*4 * \circrad, \ytop+\circrad) circle (\circrad);
    
    \draw[gray] ( \i*4 * \circrad, \ybot-\circrad) circle (\circrad);
    \draw[gray] (-\i*4 * \circrad, \ybot-\circrad) circle (\circrad);
}

\def\a{0.5}
\def\alen{0.5}
\foreach \i in {1,...,12} {
    \draw[<-, line width=1.5pt      ] ($(hatchTmid) + (\i*\a, 0)$) -- ++(0, \alen);
    \draw[<-, line width=1.5pt, gray] ($(hatchTmid) - (\i*\a, 0)$) -- ++(0, \alen);
    
    \draw[<-, line width=1.5pt, gray] ($(BhatchBmid) + (\i*\a, 0)$) -- ++(0, -\alen);
    \draw[<-, line width=1.5pt, gray] ($(BhatchBmid) - (\i*\a, 0)$) -- ++(0, -\alen);
}

\begin{scope}
  \clip (0, \ybot) rectangle (\xright, \ytop+2*\circrad+5);
  \draw[<-, line width=1.5pt] (hatchTmid) -- ++(0, \alen);
\end{scope}
\begin{scope}
  \clip (\xleft, \ybot) rectangle (0, \ytop+2*\circrad+5);
  \draw[<-, line width=1.5pt, gray] (hatchTmid) -- ++(0, \alen);
\end{scope}

\draw[<-, line width=1.5pt, gray] (BhatchBmid) -- ++(0, -\alen);

\draw ($(hatchTmid) + (0, \alen)$) -- ++(\xright+\circrad, 0) node[midway, above] {\large $F/W$};
\draw[gray] ($(hatchTmid) + (0, \alen)$) -- ++(-\xright-\circrad, 0) node[midway, above] {\large $F/W$};

\draw[gray] ($(BhatchBL) + (0, -\alen)$) -- ++(-\xleft+\xright+2*\circrad, 0) node[midway, below] {\large $F/W$};

\draw[<->, line width=1.5pt] (0, 0.5*\ybot) -- ++(\xright, 0) 
        node[midway, above] {\large W};
\draw[<->, line width=1.5pt] (0.5*\xleft, 0) -- ++(0, \ytop) 
        node[midway, right] {\large H};

\coordinate (O) at (1, 0.25);
\draw[->] (O) -- ++(1, 0) node[above] {\Large $\vec{e}_x$};
\draw[->] (O) -- ++(0, 1) node[right] {\Large $\vec{e}_y$};

\node at (\xright - 0.5, 0.5) {\Large $\Omega$};

\end{tikzpicture}
    }
    \caption{Schematic of Mandel's problem. A porous medium is vertically loaded by two plates, and is free to expand along the $x$-direction. We discretize only one quarter of the domain, exploiting symmetry.}
    \label{fig:2Dcase}
\end{figure}

We now extend our analysis to the multi-dimensional case. We choose Mandel's problem \cite{mandelConsolidationSolsEtude1953} as a two-dimensional benchmark, see \autoref{fig:2Dcase}. A rectangular domain is subjected to a consolidating force in the $y$-direction. This force is applied through rigid, impermeable plates, such that both normal displacements are constant and equal to each other along the top- and bottom edge. Along both sides, the domain can freely expand and expel pore fluid. There are two symmetry lines in this problem, allowing us to discretize only a quarter of the domain. \\
The timescales characteristic to this problem are long, in the sense that the solution is dominated by resistance to seepage, while inertia is negligible. While this makes it a suboptimal validation case from the spectral analysis perspective, it is one of the few meaningful 2D poroelasticity problems for which an analytical solution exists \cite{abousleimanMandelsProblemRevisited1996}. This analytical solution is discussed in \ref{app:mandelanalytical}.

\subsection{Validation}

\begin{figure}[h!]
    \centering
    \tikzsetnextfilename{2D_convergence}
\begin{tikzpicture}
\begin{axis}[
    width=0.8\textwidth,
    height=0.55\textwidth,
    xlabel={Number of elements along both dimensions},
    ylabel={$||U^h-U||_{L^2}$},
    xmode=log,
    ymode=log,
    legend pos=north east,
    cycle list name=mytab10,
    legend style={
        rounded corners=2pt,
        draw=black!35,
    },
    xtick pos=left,
    ytick pos=left,
    xtickten={0,1,2,3,4,5},
    minor x tick num=9,
]

\addplot table[
    x=element_count,
    y=linear_fem,
    col sep=comma
] {data/convergence2D.csv};
\addlegendentry{Linear FEM}

\addplot table[
    x=element_count,
    y=quadratic_fem,
    col sep=comma
] {data/convergence2D.csv};
\addlegendentry{Quadratic FEM}

\addplot table[
    x=element_count,
    y=cubic_fem,
    col sep=comma
] {data/convergence2D.csv};
\addlegendentry{Cubic FEM}

\addplot table[
    x=element_count,
    y=quadratic_iga,
    col sep=comma
] {data/convergence2D.csv};
\addlegendentry{Quadratic IGA}

\addplot table[
    x=element_count,
    y=cubic_iga,
    col sep=comma
] {data/convergence2D.csv};
\addlegendentry{Cubic IGA}

\end{axis}
\end{tikzpicture}
    \caption{Convergence in the $L^2$-norm of the error in fluid displacement at $t=200$ s for various discretizations of Mandel's problem, using the $u$-$U$ formulation.}
    \label{fig:convergence_mandel}
\end{figure}
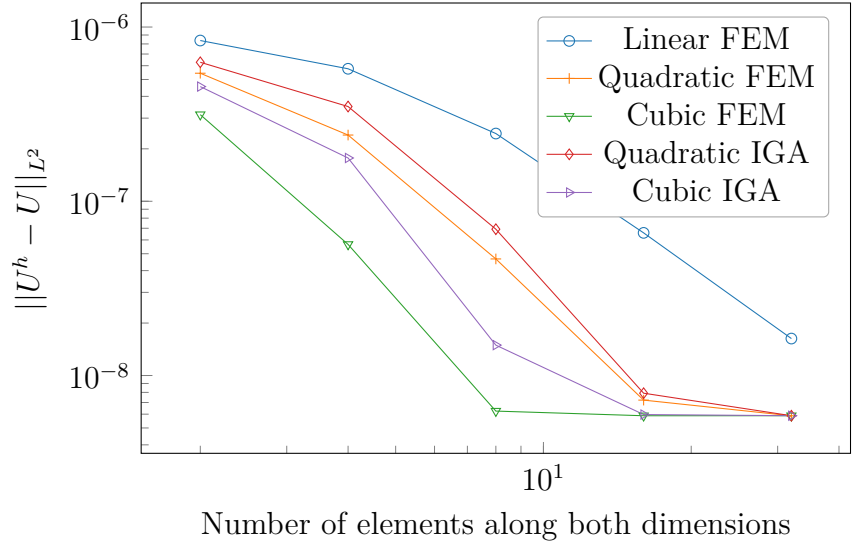

The material parameters we selected for our 2D study are listed in \autoref{tab:matpars_2D}. The magnitudes in this set approximately correspond to the anisotropic material parameters used in Ref.~\cite{abousleimanMandelsProblemRevisited1996}. Since inertia is negligible in this problem, Ref.~\cite{abousleimanMandelsProblemRevisited1996} does not define densities. Hence, we introduce arbitrary densities, chosen sufficiently low so as not to introduce inertial effects. \\
The elastic coefficients can be converted to Lamé parameters using standard relations. Additionally, $B$ and $\nu_u$ are required for the analytical solution. They are obtained through relations provided in e.g. Ref. \cite{chengPoroelasticity2016}. We choose a domain with sizes $W=$\num{0.2} \unit{\meter} and $H =$ \num{0.1} \unit{\meter}, loaded with $F=$ \num{1e7} \unit[per-mode=symbol]{\newton\per\meter}.\\
To validate our implementation, we solve Mandel's problem over 200 seconds. Time is discretized into one hundred equidistant timesteps through the implicit Newmark method with $\beta=0.25$ and $\gamma=0.5$. At the final timestep, we compute the $L^2$-norm of the difference between the numerical- and analytical fluid displacement. The results for various discretizations in space are is plotted in \autoref{fig:convergence_mandel}. In all cases, our numerical solutions converge toward the analytical one at a rate depending on the discretization until they plateau, due to limited accuracy in the time integration scheme as well as numerical precision. We consider a study of the optimality of the observed asymptotic rates beyond the scope of the current work.

{\renewcommand{\arraystretch}{1.5}
\begin{table}[h!]
    \centering
    \begin{tabular}{c|c|c|c|c|c|c|c}
        Parameter & 
        $E$ [Pa] & 
        $\nu$ [-] & 
        $ \alpha [-] $ & 
        $ \phi [-] $ & 
        $ k $ & 
        $ \rho_s [\frac{\mathrm{kg}}{\mathrm{m}^3}]$ & 
        $ \rho_f [\frac{\mathrm{kg}}{\mathrm{m}^3}]$
        \\\hline
        Value     & 
        \num{20e9} &
        \num{0.25}  & 
        \num{0.723} & 
        \num{0.02}  & 
        \num{2.0e-17} & 
        \num{1e-3} & 
        \num{1e-3}
    \end{tabular}
    \caption{Material parameters for the 2D test case. With the undrained bulk modulus $K_u = K + \alpha^2 Q$, Skempton's coefficient is $B = \alpha Q/K_u$ and the undrained Poisson ratio is $\nu_u = \frac{3K_u-2\mu}{2(3K_u+\mu)}$.}
    \label{tab:matpars_2D}
\end{table}}

\subsection{Spectral analysis}
\begin{figure}[h!]
    \centering
    \def\markersize{3.0pt}
\def\markerlinewidth{1.0pt}

\tikzsetnextfilename{2D_spectra_upU_vs_uU}
\begin{tikzpicture}
\begin{groupplot}[
    group style={
        group size=2 by 1,
        horizontal sep=0.06\textwidth,
        y descriptions at=edge left,
    },
    width=0.495\textwidth,
    height=0.495\textwidth,
    xlabel={$n$ [-]},
    ylabel={$\omega_n / \omega_{n, \mathrm{ref}}$},
    ymode=log,
    ymax=20,
    xtick pos=left,
    ytick pos=left,
    cycle list name=mytab10,
    legend image post style={
        mark indices={},
    },
    legend style={
        rounded corners=2pt,
        draw=black!35,
        cells={anchor=west},
    },
]

\nextgroupplot[
    legend to name=leftlegend,
    legend columns=2,
    xtick={0,300,600,900,1200},    
]

\addplot+[
    mark=o,
    mark indices={0,25,...,10000},
    mark size=\markersize,
    mark options={
        fill=none,
        line width=\markerlinewidth,
        solid,
    },
] table[
    x=mode_number,
    y=eigfreq_over_reference,
    col sep=comma
] {data/2D_spectra_divfiltered/2D_spectrum_cubic_uUconforming_C2_upU_divfiltered_normalized.csv};
\addlegendentry{$\alpha = 2$, $u$-$p$-$U$}

\addplot+[
    mark=x,
    mark indices={0,25,...,10000},
    mark size=\markersize,
    mark options={
        fill=none,
        line width=\markerlinewidth,
        solid,
    },
] table[
    x=mode_number,
    y=eigfreq_over_reference,
    col sep=comma
] {data/2D_spectra_divfiltered/2D_spectrum_cubic_uUconforming_C2_uU_divfiltered_normalized.csv};
\addlegendentry{$\alpha = 2$, $u$-$U$}

\addplot+[
    mark=o,
    mark indices={0,50,...,10000},
    mark size=\markersize,
    mark options={
        fill=none,
        line width=\markerlinewidth,
        solid,
    },
] table[
    x=mode_number,
    y=eigfreq_over_reference,
    col sep=comma
] {data/2D_spectra_divfiltered/2D_spectrum_cubic_uUconforming_C1_upU_divfiltered_normalized.csv};
\addlegendentry{$\alpha = 1$, $u$-$p$-$U$}

\addplot+[
    mark=x,
    mark indices={0,50,...,10000},
    mark size=\markersize,
    mark options={
        fill=none,
        line width=\markerlinewidth,
        solid,
    },
] table[
    x=mode_number,
    y=eigfreq_over_reference,
    col sep=comma
] {data/2D_spectra_divfiltered/2D_spectrum_cubic_uUconforming_C1_uU_divfiltered_normalized.csv};
\addlegendentry{$\alpha = 1$, $u$-$U$}

\nextgroupplot[
    legend to name=rightlegend,
    legend columns=1,
]

\addplot+[
    mark=o,
    mark indices={0,25,...,10000},
    mark size=\markersize,
    mark options={
        fill=none,
        line width=\markerlinewidth,
        solid,
    },
] table[
    x=mode_number,
    y=eigfreq_over_reference,
    col sep=comma
] {data/2D_spectra_divfiltered/2D_spectrum_cubic_Uconforming_uiso_C2_upU_divfiltered_normalized.csv};
\addlegendentry{$u$-$p$-$U$}

\addplot+[
    mark=x,
    mark indices={0,25,...,10000},
    mark size=\markersize,
    mark options={
        fill=none,
        line width=\markerlinewidth,
        solid,
    },
] table[
    x=mode_number,
    y=eigfreq_over_reference,
    col sep=comma
] {data/2D_spectra_divfiltered/2D_spectrum_cubic_Uconforming_uiso_C2_uU_divfiltered_normalized.csv};
\addlegendentry{$u$-$U$}

\end{groupplot}

\node[anchor=north] at ($(group c1r1.south)+(0,-1.2cm)$)
    {\pgfplotslegendfromname{leftlegend}};

\node[anchor=north] at ($(group c2r1.south)+(0,-1.2cm)$)
    {\pgfplotslegendfromname{rightlegend}};

\node[anchor=north] at ($(group c1r1.south)+(0,-2.8cm)$)
    {Divergence-conforming $\trialspaceu$};

\node[anchor=north] at ($(group c2r1.south)+(0,-2.8cm)$)
    {Isotropic $\trialspaceu$};

\end{tikzpicture}
    \caption{Normalized spectra for the two formulations, not including divergence free modes. In (a)~both displacement fields are discretized with divergence-conforming elements at varying regularity. In (b)~the matrix displacement is discretized using isotropic elements instead, and only optimal regularity elements are used.}
    \label{fig:2D_upU_vs_uU}
\end{figure}
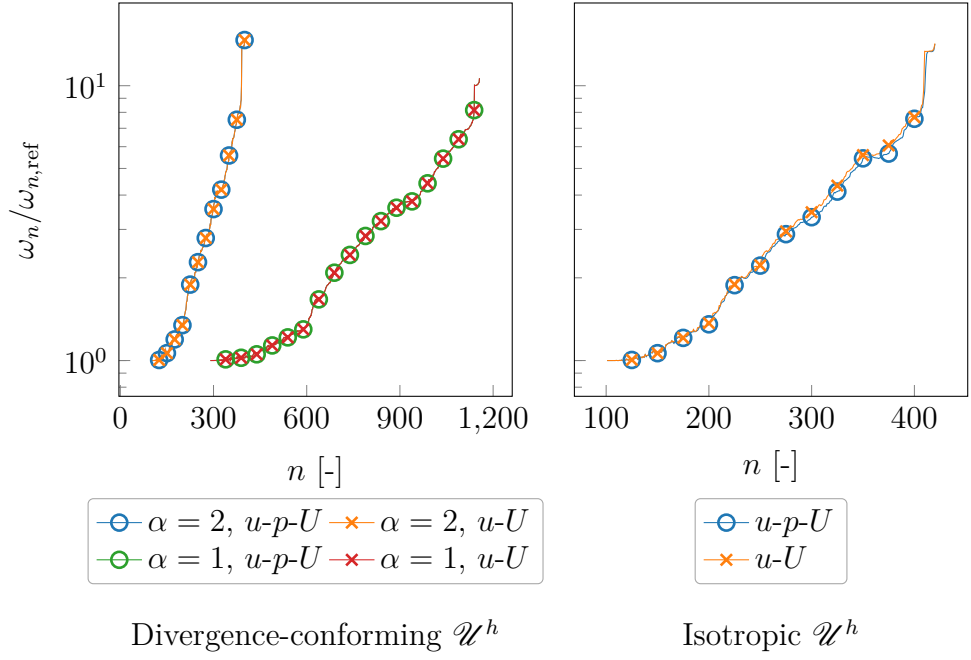

To study the discrete spectra in 2D, we use the same geometry as that of the Mandel problem. The same symmetry boundaries are used as well, but we no longer linearly constrain the normal displacements along the top edge to be equal, such that the boundary conditions are separable.\\
With that, we first study spectra of the $u$-$U$ formulation. In \autoref{fig:eigs2D_nonconforming}, we compare two spectra resulting from cubic, optimal regularity elements. In the one case, the fluid displacement is discretized with divergence-conforming elements. That is, according to \eqref{eq:conforming2D-vector} with $k_1=k_2=3$ and $\alpha_1=\alpha_2 = 2$. In the other case, we use isotropic (bi-cubic) elements. In both cases, isotropic elements are used to discretize the matrix displacement. 
For the conforming discretization, we see a distinct null space. The eigenmodes spanning this space describe divergence-free fluid motions. In the strong form \eqref{eq:twofield_strong} it is apparent that such modes are restricted only through damping. Since we here study the undamped problem, they form a null space. We did confirm that a perturbation by the damping matrix makes the problem positive definite. The $L^2$-norm of the divergence of $\fludisp$ in each mode is pictured in \autoref{fig:2D_div_per_mode}. We see one hundred modes with a negligible divergence. This is exactly the expected count: the divergence of the trial space is $\mathbb{S}^{2}_{1} \otimes \mathbb{S}^2_1$, which for $8\times8$ elements contains $10\times10$ DOFs. Conversely, for the isotropic discretization there is no well-separated null space. Only nine eigenmodes are approximately zero. We conclude that, to correctly capture the structure inherent in the weak form, one should choose a divergence conforming discretization for the fluid displacement. Note that it is not necessary to use one for $\matdisp$: whether $\trialspaceu$ is isotropic or conforming does not influence the size of the null space.\\
\begin{figure}[h!]
    \centering
    \def\markersize{1.5pt}

\tikzsetnextfilename{2D_spectra_nonconforming}
\begin{tikzpicture}
\begin{axis}[
    width=0.9\textwidth,
    height=0.55\textwidth,
    xlabel={$n$ [-]},
    ylabel={$\sqrt{||\omega^2||}$ [\unit{\radian}/\unit{\second}]},
    ymode=log,
    legend pos=south east,
    legend style={
        rounded corners=2pt,
        draw=black!35,
        cells={anchor=west},
    },
    xtick pos=left,
    ytick pos=left,
]

\addplot+[
    only marks,
    mark=*,
    mark size=\markersize,
    mark options={
        draw=tabblue,
        fill=none,
    },
    color=tabblue,
] table[
    x=mode_number,
    y=eigfreq_norm,
    col sep=comma
] {data/2D_spectra/2D_spectrum_cubic_iso_C2_uU.csv};
\addlegendentry{Isotropic $\trialspaceu$, isotropic $\trialspaceU$}

\addplot+[
    only marks,
    mark=*,
    mark size=\markersize,
    mark options={
        draw=taborange,
        fill=none,
    },
    color=taborange,
] table[
    x=mode_number,
    y=eigfreq_norm,
    col sep=comma
] {data/2D_spectra/2D_spectrum_cubic_Uconforming_C2_uU.csv};
\addlegendentry{Isotropic $\trialspaceu$, div-conforming $\trialspaceU$}

\end{axis}
\end{tikzpicture}
    \caption{Eigenfrequencies resulting from optimal regularity cubic discretizations of the $u$-$U$ formulation. Obtained with an isotropic discretization for $\matdisp$ and either an isotropic or a divergence-conforming discretization for $\fludisp$. Discretized with $8\times8$ elements.}
    \label{fig:eigs2D_nonconforming}
\end{figure}
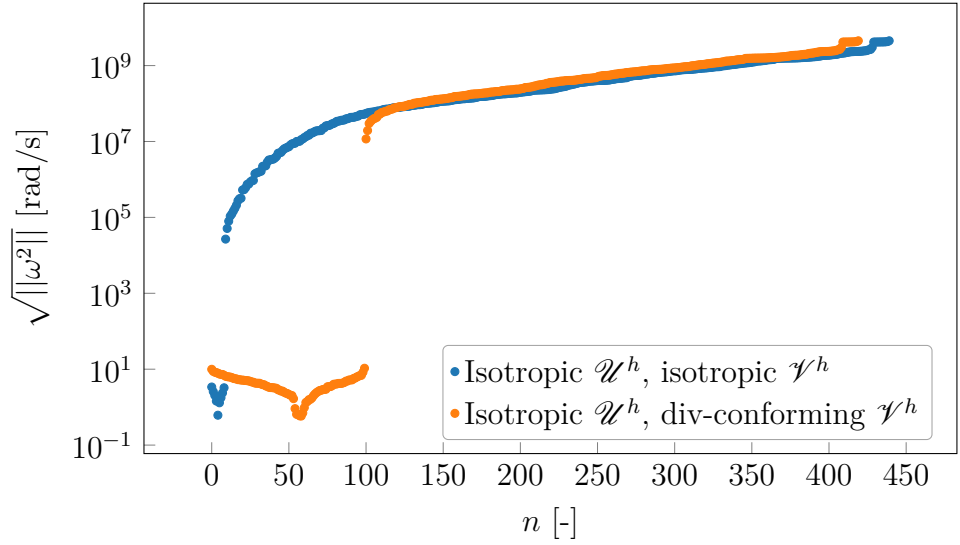
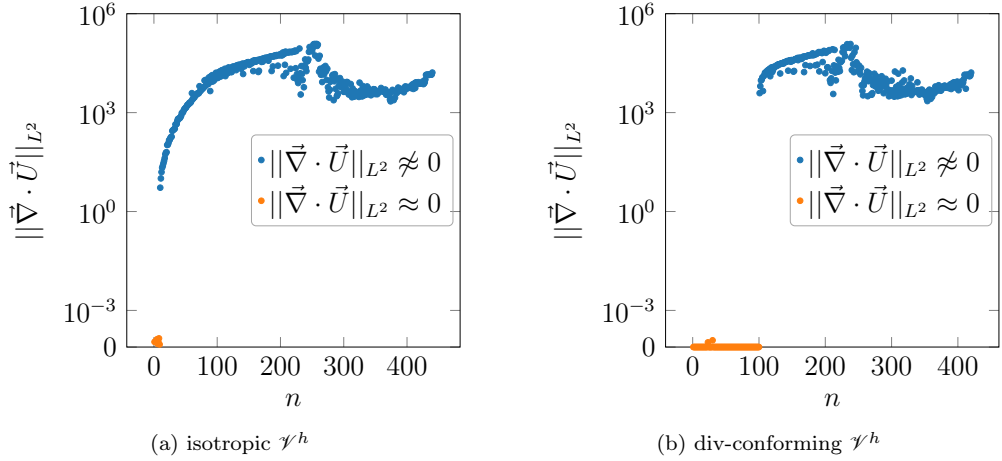
\begin{figure}[h!]
    \centering
     \begin{subfigure}[t]{0.48\textwidth}
        \centering
        \resizebox{\linewidth}{!}{%
            \def\markersize{1.2pt}

\tikzsetnextfilename{2D_div_per_mode_iso}
\begin{tikzpicture}
\begin{axis}[
  width=1.0\linewidth,
  height=1.0\linewidth,
  xlabel={$n$},
  ylabel={$||\grad \cdot \fludisp||_{L^2}$},
  legend style={
    at={(0.98,0.50)},
    anchor=east,
    fill=none,
    rounded corners=2pt,
    draw=black!35,
  },
  ymode=normal,
  ymin=0.0,
  ymax= 1.011111111111111e-02,
  scaled y ticks=false,
  ytick={
    0.0,
    1.111111111111111e-03,
    4.111111111111111e-03,
    7.111111111111111e-03,
    1.011111111111111e-02},
  yticklabels={
    $0$,
    $10^{-3}$,
    $10^{0}$,
    $10^{3}$,
    $10^{6}$,
  }
]

\addplot[
  scatter,
  only marks,
  scatter src=explicit symbolic,
  mark=*,
  mark size=\markersize,
  scatter/classes={
    0={tabblue, mark options={solid, fill=tabblue}},
    1={taborange, mark options={solid, fill=taborange}}
  }
]
table[
  x=mode_number,
  y=divU_symlog,
  meta=is_approximately_zero_divU,
  col sep=comma
]{data/2D_div_per_mode/2D_spectrum_cubic_iso_C2_divU.csv};

\addlegendentry{$||\grad\cdot\fludisp||_{L^2} \not\approx 0$}
\addlegendentry{$||\grad\cdot\fludisp||_{L^2} \approx 0$}

\end{axis}
\end{tikzpicture}
        }
        \caption{isotropic $\trialspaceU$}
        \label{fig:2D_div_per_mode_a}
    \end{subfigure}
    \hfill
    \begin{subfigure}[t]{0.48\textwidth}
        \centering
        \resizebox{\linewidth}{!}{%
            \def\markersize{1.2pt}

\tikzsetnextfilename{2D_div_per_mode_conforming}
\begin{tikzpicture}
\begin{axis}[
  width=1.0\linewidth,
  height=1.0\linewidth,
  xlabel={$n$},
  ylabel={$||\grad \cdot \fludisp||_{L^2}$},
  legend style={
    at={(0.98,0.50)},
    anchor=east,
    fill=none,
    rounded corners=2pt,
    draw=black!35,
  },
  ymode=normal,
  ymin=0.0,
  ymax= 1.011111111111111e-02,
  scaled y ticks=false,
  ytick={
    0.0,
    1.111111111111111e-03,
    4.111111111111111e-03,
    7.111111111111111e-03,
    1.011111111111111e-02},
  yticklabels={
    $0$,
    $10^{-3}$,
    $10^{0}$,
    $10^{3}$,
    $10^{6}$,
  }
]

\addplot[
  scatter,
  only marks,
  scatter src=explicit symbolic,
  mark=*,
  mark size=\markersize,
  scatter/classes={
    0={tabblue, mark options={solid, fill=tabblue}},
    1={taborange, mark options={solid, fill=taborange}}
  }
]
table[
  x=mode_number,
  y=divU_symlog,
  meta=is_approximately_zero_divU,
  col sep=comma
]{data/2D_div_per_mode/2D_spectrum_cubic_Uconforming_C2_divU.csv};    

\addlegendentry{$||\grad\cdot\fludisp||_{L^2} \not\approx 0$}
\addlegendentry{$||\grad\cdot\fludisp||_{L^2} \approx 0$}

\end{axis}
\end{tikzpicture}
        }
        \caption{div-conforming $\trialspaceU$}
        \label{fig:2D_div_per_mode_b}
    \end{subfigure}
    \caption{The $L^2$-norm of the divergence of $\fludisp$, for each mode in \autoref{fig:eigs2D_nonconforming}.}
    \label{fig:2D_div_per_mode}
\end{figure}
We proceed by comparing spectra for the two formulations in \autoref{fig:2D_upU_vs_uU}. Since we do not have an analytical solution, normalization is not straightforward. We choose to plot only non-zero eigenmodes, and normalize those against the non-zero eigenmodes in a refined discretization. This refined case consists of $32\times32$ fifth-order, optimal regularity, divergence-conforming elements for both fields. The other cases use $8\times8$ cubic elements. \autoref{fig:2D_upU_vs_uU}a shows the cases wherein both fields are discretized with div-conforming spaces, while \autoref{fig:2D_upU_vs_uU}b shows cases where the matrix displacement is discretized by isotropic elements. Together, these confirm, now in 2D, our expectations set in \autoref{subsec:conforming}: the two formulations coincide if and only if the divergence of any function in $\trialspaceU$ or $\trialspaceu$ is in $\trialspacep$. \\
Since the formulations perform identically for compressible problems with compatible spaces, and the $u$-$p$-$U$ formulation requires additional computational effort in inverting the matrix $\matr{P}$, we conclude that the $u$-$U$ formulation is preferrable in this context.

\begin{remark}
    Note that the normalization in \autoref{fig:2D_upU_vs_uU} does not provide a measure of discretization error. We attribute this to the presence of shear modes and longitudinal modes, which we cannot cleanly separate. Consider the case where shear modes are typically lower in eigenfrequency than longitudinal modes. When we then compare all $n$ modes in a given discretization against the lowest $n$ modes in the refined reference, the latter set could consist of exclusively shear modes. Therefore, we cannot expect this normalization to result only in values close to one.
\end{remark}

\subsection{The critical timestep}
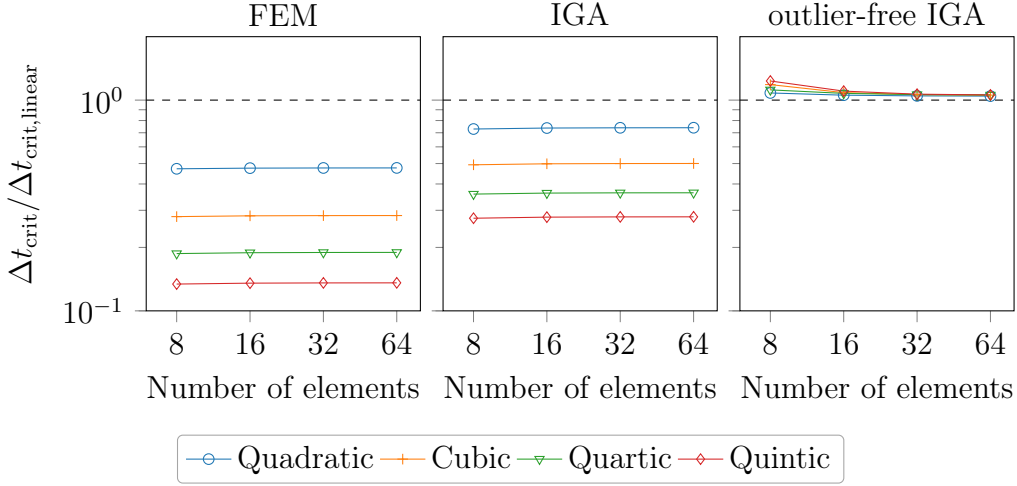
\begin{figure}[h!]
    \centering

    \def\marksize{2pt}
    \tikzsetnextfilename{2D_critstepscaling}
\begin{tikzpicture}

\begin{groupplot}[
    group style={
        group size=3 by 1,
        horizontal sep=0.3cm,
    },
    width=.38\textwidth,
    height=.38\textwidth,
    grid=none,
    xmode=log,
    ymode=log,
    xmin=6,
    xmax=80,
    ymin=0.1,
    ymax=2.0,
    xlabel={Number of elements},
    cycle list name=mytab10,
    tick align=outside,
    tick pos=left,
    xtick={8, 16, 32, 64},
    xticklabels={$8$,$16$,$32$,$64$},
]

\nextgroupplot[
    title={FEM},
    ylabel={$\Delta t_\mathrm{crit} / \Delta t_\mathrm{crit,linear}$},
    legend columns=5,
    legend to name=2D_critstepscalinglegend,
    legend image post style={
        mark indices={2}
    },
    legend style={
        rounded corners=2pt,
        draw=black!35,
    },
]

\addplot[
    dashed,
    black,
    forget plot,
] coordinates {(6,1.0) (80,1.0)};

\addplot+[
    mark size=\marksize,
]
table[col sep=comma, x=element_count, y=quadratic_normalized_by_linear]
{data/2D_maxeigfreqinverse/2D_maxeigfreqinverse_div_conforming__fem.csv};
\addlegendentry{Quadratic}

\addplot+[
    mark size=\marksize,
]
table[col sep=comma, x=element_count, y=cubic_normalized_by_linear]
{data/2D_maxeigfreqinverse/2D_maxeigfreqinverse_div_conforming__fem.csv};
\addlegendentry{Cubic}

\addplot+[
    mark size=\marksize,
]
table[col sep=comma, x=element_count, y=quartic_normalized_by_linear]
{data/2D_maxeigfreqinverse/2D_maxeigfreqinverse_div_conforming__fem.csv};
\addlegendentry{Quartic}

\addplot+[
    mark size=\marksize,
]
table[col sep=comma, x=element_count, y=quintic_normalized_by_linear]
{data/2D_maxeigfreqinverse/2D_maxeigfreqinverse_div_conforming__fem.csv};
\addlegendentry{Quintic}

\nextgroupplot[
    title={IGA},
    yticklabels={},
]

\addplot[
    dashed,
    black,
    forget plot,
] coordinates {(6,1.0) (80,1.0)};

\addplot+[
    mark size=\marksize,
]
table[col sep=comma, x=element_count, y=quadratic_normalized_by_linear]
{data/2D_maxeigfreqinverse/2D_maxeigfreqinverse_div_conforming__iga.csv};

\addplot+[
    mark size=\marksize,
]
table[col sep=comma, x=element_count, y=cubic_normalized_by_linear]
{data/2D_maxeigfreqinverse/2D_maxeigfreqinverse_div_conforming__iga.csv};

\addplot+[
    mark size=\marksize,
]
table[col sep=comma, x=element_count, y=quartic_normalized_by_linear]
{data/2D_maxeigfreqinverse/2D_maxeigfreqinverse_div_conforming__iga.csv};

\addplot+[
    mark size=\marksize,
]
table[col sep=comma, x=element_count, y=quintic_normalized_by_linear]
{data/2D_maxeigfreqinverse/2D_maxeigfreqinverse_div_conforming__iga.csv};

\nextgroupplot[
    title={outlier-free IGA},
    yticklabels={},
]

\addplot[
    dashed,
    black,
    forget plot,
] coordinates {(6,1.0) (80,1.0)};

\addplot+[
    mark size=\marksize,
]
table[col sep=comma, x=element_count, y=quadratic_normalized_by_linear]
{data/2D_maxeigfreqinverse/2D_maxeigfreqinverse_div_conforming__outlierfree_iga.csv};

\addplot+[
    mark size=\marksize,
]
table[col sep=comma, x=element_count, y=cubic_normalized_by_linear]
{data/2D_maxeigfreqinverse/2D_maxeigfreqinverse_div_conforming__outlierfree_iga.csv};

\addplot+[
    mark size=\marksize,
]
table[col sep=comma, x=element_count, y=quartic_normalized_by_linear]
{data/2D_maxeigfreqinverse/2D_maxeigfreqinverse_div_conforming__outlierfree_iga.csv};

\addplot+[
    mark size=\marksize,
]
table[col sep=comma, x=element_count, y=quintic_normalized_by_linear]
{data/2D_maxeigfreqinverse/2D_maxeigfreqinverse_div_conforming__outlierfree_iga.csv};

\end{groupplot}

\end{tikzpicture}
    \vspace{-0.5em}

    \begingroup
    \tikzexternaldisable
    \pgfplotslegendfromname{2D_critstepscalinglegend}
    \endgroup

    \vspace{0.5em}
    \caption{In 2D, critical timesteps relative to those for linear FEM. Using the $u$-$U$ formulation, discretized with FEM, IGA, and outlier-free IGA at various orders.}
    \label{fig:2D_critstepscaling}
\end{figure}

In \autoref{fig:2D_critstepscaling}, we plot how the critical timestep scales with order for various discretizations. Again, this mirrors the findings for elasticity in Ref. \cite{hiemstraRemovalSpuriousOutlier2021}: the critical timestep scales strongly with the polynomial order in FEM. This scaling is less strong with IGA, and virtually eliminated by outlier-free IGA. Note that \autoref{fig:2D_critstepscaling} was created with divergence-conforming spaces. For the $u$-$U$ formulation, using an isotropic discretization either only for $\matdisp$ or for both displacements gives virtually equal results with regards to scaling of the timestep.

\section{Conclusion}\label{chap:conc}

We investigated whether the advantageous properties of higher-order IGA as compared to higher-order FEM, which have been well-established in elastodynamics, extend to poroelastodynamics. Through studying eigenfrequency spectra of poroelasticity for various discretizations, we have found that this is indeed the case. \\
In 1D, when we split the spectrum into fast and slow waves, the spectra carry the exact same characteristics as those resulting from elastodynamics. Spectra resulting from higher-order FEM include non-physical optical modes, while higher-order IGA computes only acoustic modes. That is, with the exception of a small number of outlier modes. These can be readily eliminated through methods established in the literature, which effectively results in the usage of optimal spline spaces to approximate the physical fields. \\
We showed that the critical timestep in an explicit solver depends on the polynomial order. When using FEM, the critical timestep decreases strongly upon order elevation. This effect is reduced using IGA, and essentially eliminated with outlier removal. Thus, optimal spline spaces are highly promising for explicit dynamics in terms of computational efficiency, as they eliminate one of the costs usually associated with higher-order methods.\\
We also compared two poroelasticity formulations in 1D. We proved and numerically demonstrated that, provided that divergence-conforming approximation spaces are chosen, both formulations result in identical spectra. In this context, we recommend the $u$-$U$ formulation, as it is computationally less expensive than the $u$-$p$-$U$ formulation. Overall, we can conclude that IGA, especially with outlier elimination, is promising as an efficient method for dynamic poroelasticity. It provides higher-order discretizations without inflating the number of DOFs, or decreasing the admissible timestep. \\
We note that higher-order methods are favorable only when the solution is sufficiently smooth. This is typically the case in acoustics, but heterogeneity could be limiting in our application. In addition, IGA can only compete with FEM in terms of computational effort when it is implemented efficiently, and so long as the effort involved in numerical integration does not become prohibitive. We consider an efficient implementation to be a key next step in employing IGA and poroelasticity for seismic hazard assessment at scale. In addition, we intend to study the accuracy with which poroelasticity captures near-surface ground motion.

\section*{Acknowledgments} 
This work is supported by the Nederlandse Organisatie voor Wetenschappelijk Onderzoek (NWO) DeepNL program, Grant DEEP.NL.2023.018. René Hiemstra acknowledges funding from the European Union's Horizon 2020 research and innovation programme under the Marie Skłodowska-Curie grant agreement No 101105786.

\newpage
\appendix

\section{Requirements for equivalent formulations}\label{app:equivalence}
The $u$-$p$-$U$ formulation, after eliminating the pressure, only differs from the $u$-$U$ formulation by the stiffness matrix (see Eq. \eqref{eq:secondorderform}). That is, provided that for both formulations the same discretization is chosen for the fields $\matdisp$ and $\fludisp$. Then, the formulations coincide in the case that 
\begin{equation}
    \overline{\matr{K}}_{u\mathrm{-}p\mathrm{-}U} =     \overline{\matr{K}}_{u\mathrm{-}U}.
\end{equation}
Substituting their definitions \eqref{eq:upUstiff}, \eqref{eq:uUstiff}, these matrices are equal when 
\begin{subequations}
\begin{align}
    \label{eq:equalityreq_step1_eq1}
    \matr{G}\matr{P}^{-1}\matr{G^T} = \matr{R}_1,    \\
    \label{eq:equalityreq_step1_eq2}
    \matr{G}\matr{P}^{-1}\matr{H^T} = \matr{R}_2,    \\
    \label{eq:equalityreq_step1_eq3}
    \matr{H}\matr{P}^{-1}\matr{H^T} = \matr{R}_3.
\end{align}
\end{subequations}
We first consider the condition \eqref{eq:equalityreq_step1_eq1}. Upon substituting the definitions of the operators \eqref{eq:discreteoperators}, \eqref{eq:operators1} and \eqref{eq:operators2}, we obtain the requirement
\begin{align}\label{eq:equalityreq_step2}
    (\Pi_p d_u, 
     \Pi_p e_u)_{L^2} = 
    (d_u, 
     e_u)_{L^2} \quad \forall\, d_u, e_u \in \grad \cdot \trialspaceu,
\end{align}
where $\Pi_p$ is the $L^2$ projection onto $\trialspacep$. Any divergence is equal to the sum of its projection and the part missed by the projection, \emph{i.e.}
\begin{equation}
    d_u = \Pi_p d_u + (\mathcal{I} - \Pi_p) d_u \quad \forall \, d_u \in \grad \cdot \trialspaceu,
\end{equation}
with $\mathcal{I}$ the identity operator. By orthogonality of the $L^2$ projection, we have that 
\begin{align}
    ||                      d_u ||_{L^2}^2 =
    ||               \Pi_p  d_u ||_{L^2}^2 +
    ||(\mathcal{I} - \Pi_p) d_u ||_{L^2}^2 \quad \forall \, d_u \in \grad \cdot \trialspaceu.
\end{align}
Combined with \eqref{eq:equalityreq_step2}, it follows that \eqref{eq:equalityreq_step1_eq1} is fulfilled if and only if the part missed by the projection vanishes, or equivalently
\begin{equation}
    d_u = \Pi_p d_u \quad \forall \, d_u \in \grad \cdot \trialspaceu.
\end{equation}
That is, the divergence of $\trialspaceu$ must be exactly represented by its $L^2$ projection onto $\trialspacep$, \emph{i.e.}
\begin{equation}
    \grad \cdot \trialspaceu \subseteq \trialspacep.
\end{equation}
The third requirement, \eqref{eq:equalityreq_step1_eq3}, is isomorphic to the first. Hence, it is fulfilled when the divergence of $\trialspaceU$ lies within $\trialspacep$. If we already pose both of these requirements, the second, mixed statement \eqref{eq:equalityreq_step1_eq2} follows.

\section{Damped spectra}\label{app:damped}

\begin{figure}[h!]
    \centering
    
    \def\marksize{1pt}
    \tikzsetnextfilename{complex_plane}
\begin{tikzpicture}

\begin{axis}[
    width=.8\textwidth,
    height=.5\textwidth,
    grid=both,
    xlabel={$Re(\omega)$},
    ylabel={$Im(\omega)$},
    ymin=-17.118095651,
    ymax= 17.118095651,
    xmin=-70,
    xmax=5,
    enlargelimits=false,
    scaled y ticks=false,
    ytick={
        -12.512925465,
        -7.907755279,
        -3.302585093,
         0,
         3.302585093,
         7.907755279,
         12.512925465
    },
    yticklabels={
        {$-10^7$},
        {$-10^5$},
        {$-10^3$},
        {$0$},
        {$10^3$},
        {$10^5$},
        {$10^7$}
    },
]

\addplot[
    only marks,
    mark size=\marksize,
    color=tabblue,
]
table[
    col sep=comma,
    x=real_omega,
    y=symlog_imag_omega
]
{data/damped_p5-5_c4-4.csv};

\end{axis}

\end{tikzpicture}
        
    \vspace{0.25em}
    \caption{Complex eigenfrequencies resulting from the damped poroelasticity system. Obtained using the $u$-$U$ formulation with fifth order optimal regularity splines, using 512 elements.}
    \label{fig:complexfrequencies}
\end{figure}
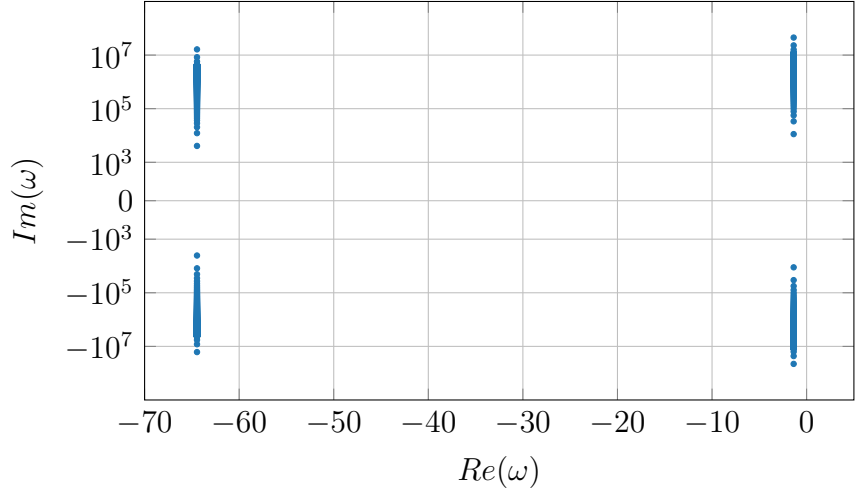

For the spectral analysis in \cref{chap:1Dcase,chap:2Dcase}, we neglected damping to obtain a generalized eigenvalue problem. When we include damping, this is no longer the case. Instead, after substituting the ansatz $[\matr{u}, \matr{U}]^T = \matr{x} \exp(\omega t)$ into \eqref{eq:secondorderform}, we obtain the quadratic eigenvalue problem

\begin{equation}\label{eq:QEP}
    \left[\omega^2 \overline{\matr{M}} + \omega \overline{\matr{C}} + \overline{\matr{K}}\right]\matr{x} = 0.
\end{equation}
For an overview of solution strategies, we refer to Ref. \cite{tisseurQuadraticEigenvalueProblem2001}. The general approach is to linearize the system, resulting in a generalized eigenvalue problem twice as large as the original. Ref. \cite{tisseurQuadraticEigenvalueProblem2001} describes two general linearizations, the first and second companion form. While it is possible to choose these such that the system remains symmetric, in our underdamped case, it will not be positive definite. Since we are not aware of eigensolvers that exploit symmetry without requiring positive definiteness, we do not pick a symmetric linearization. Rather, we choose it such that the right-hand side matrix becomes identity, thereby reducing further to a standard eigenvalue problem. To that end, we first eliminate the mass matrix through the Cholesky decomposition and a substitution,

\begin{equation}\label{eq:llt}
    \overline{\matr{M}} = \matr{L}\matr{L}^T, \qquad \matr{x} = \matr{L}^{-T}\matr{y}.
\end{equation}
If we left-multiply \autoref{eq:QEP} by $\matr{L}^{-1}$ and insert our substition, we obtain
\begin{equation}\label{eq:QEP_reduced}
    \left[\omega^2 \matr{I} + \omega \hat{\matr{C}} + \hat{\matr{K}}\right]\matr{y} = 0,
\end{equation}
where $\hat{\matr{C}} = \matr{L}^{-1}\overline{\matr{C}}\matr{L}^{-T}$ and $\hat{\matr{K}} = \matr{L}^{-1}\overline{\matr{K}}\matr{L}^{-T}$. Then, linearizing by the first companion form yields
\begin{equation} \label{eq:linearizedQEP}
    \begin{bmatrix}
        \matr{0} && \matr{I} \\
        -\hat{\matr{K}} && -\hat{\matr{C}}
    \end{bmatrix}
    \begin{bmatrix}
        \matr{y} \\
        \omega \matr{y}
    \end{bmatrix}
    -\omega
    \begin{bmatrix}
        \matr{I} && \matr{0} \\
        \matr{0} && \matr{I}
    \end{bmatrix}
    \begin{bmatrix}
        \matr{y} \\
        \omega \matr{y}
    \end{bmatrix}
    =\matr{0}.
\end{equation}
Since the right-hand side is now the identity, this is a standard eigenvalue problem. Solving with a standard solver for indefinite systems yields a full set of $2n$ complex eigenvalues and vectors. The upper half of the resulting eigenvectors is independent, as is visible in \autoref{eq:linearizedQEP}. This is transformed back into the original coordinate system according to \autoref{eq:llt}. Note that if the transformed vectors were normalized such that $\matr{y}^T \matr{y} = 1$, the eigenvectors themselves are normalized in the mass norm, $\matr{x} \overline{\matr{M}} \matr{x}^T = 1$.
\\
This procedure was applied to a single discretization of the 1D eigenvalue problem studied in \autoref{chap:1Dcase}. The resulting complex eigenfequencies are provided in \autoref{fig:complexfrequencies}, where we find two coherent clusters, each with a practically constant real part. Thus, we have strongly and weakly damped branches, which we interpret to correspond to slow and fast waves respectively. We note that refining the discretization yields additional modes within these clusters, but does not shift them. \\
The imaginary part, especially that of the highest modes, is orders of magnitude larger than their real part, meaning that $|\omega| \approx |Im(\omega)|$. Consequently, the critical timestep can be predicted very well through the eigenfrequencies of the undamped problem.

\section{Analytical solution to Mandel's problem}\label{app:mandelanalytical}

The analytical solution as described by \cite{abousleimanMandelsProblemRevisited1996} is valid for transversely isotropic media. Here, we elaborate it in the fully isotropic setting. \\
The solution is a sum over values $\beta$, which satisfy the characteristic function
\begin{equation}
    \tan(\beta) = \beta \frac{1-\nu}{\nu_u - \nu},
\end{equation}
where $\nu$ is the Poisson ratio, and $\nu_u$ is the undrained Poisson ratio. Next, a generalized consolidation coefficient is defined, as
\begin{equation}
    c = 2 k B^2\mu\frac{(1-\nu)(1+\nu_u)^2}{9(1-\nu_u)(\nu_u-\nu)},
\end{equation}
where $B$ is Skempton's coefficient \cite{skemptonPorePressureCoefficients1954}. With that, the $x$-component of the displacement, which is independent of the $y$-coordinate, is given as
\begin{align}
    \begin{split}
        u_x(x, t) = \frac{F\nu} {2 \mu} - \frac{F \nu_u}{\mu} \sum_\beta \left[\frac{\sin(\beta)\cos(\beta)}{\beta-\sin(\beta)\cos(\beta)}\exp\left(\frac{\beta^2 ct}{W^2}\right)\right]x + \\
        \frac{FW}{\mu}\sum_\beta\left[ \frac{\cos(\beta)}{\beta - \cos(\beta)\sin(\beta)} \sin\left(\frac{\beta x}{W}\right)\exp\left(\frac{\beta^2 ct}{W^2}\right)   \right],
    \end{split}
\end{align}
where W and F are respectively the width of the domain and the distributed load, as indicated in \autoref{fig:2Dcase}. The $y$-component of the displacement is independent of the $x$-coordinate, and is expressed as
\begin{align}
    \begin{split}
        u_y(y, t) = -\frac{F(1-\nu)}{2\mu}& + \\
        \frac{F(1-\nu_u)}{\mu}&\sum_\beta \left[\frac{\sin(\beta)\cos(\beta)}{\beta-\sin(\beta)\cos(\beta)} \exp\left(\frac{\beta^2 ct}{W^2}\right)\right] y.
    \end{split}
\end{align}
For the fluid field, the analytical solution is given in terms of the Darcy flux. By integrating this over time, we obtain a solution for $\vec{w}$. The $y$-component of $\vec{w}$ vanishes, and the $x$-component, which is independent of the $y$-coordinate, is expressed as
\begin{align}
    \begin{split}            
        w_x(x, t) = \frac{2FBk(1+\nu_u)}{3W} \times\\\sum_\beta\left[
        \left(1 - \exp\left(\frac{\beta^2 ct}{W^2}\right) \right) \frac{W^2}{\beta^2 c} \frac{\beta\sin(\beta)}{\beta - \sin(\beta)\cos(\beta)            }\sin\left(\frac{\beta x}{W}\right)
        \right].
    \end{split}
\end{align}

\newpage

\bibliographystyle{elsarticle-num} 
\bibliography{bibliography}

\end{document}